\documentclass[11pt]{amsart}
\usepackage{amssymb,latexsym}
\input diagrams.tex
\diagramstyle[small,labelstyle=\scriptstyle,midshaft]

\theoremstyle{plain}
\newtheorem{thm}{Theorem}[section]
\newtheorem{cor}[thm]{Corollary}
\newtheorem{lem}[thm]{Lemma}
\newtheorem{pro}[thm]{Proposition}

\theoremstyle{definition}
\newtheorem{block}[thm]{}

\theoremstyle{remark}
\newtheorem{rem}{Remark}
\newtheorem*{notation}{Notation and conventions} 
\newtheorem{stepa}{Step}
\newtheorem{step}{Step}

\newcommand{\lto}{\longrightarrow}

\DeclareMathOperator{\codim}{codim}

\begin{document}
\title[Hurwitz spaces of quadruple coverings]
{\bf Hurwitz spaces of quadruple coverings of elliptic curves
 and the moduli space of abelian threefolds $\mathcal{A}_3(1,1,4)$
 }
\begin{abstract}
We prove that the moduli space $\mathcal{A}_3(1,1,4)$ of polarized abelian threefolds with polarization of type $(1,1,4)$ is unirational. By a result of Birkenhake and Lange this implies the unirationality of the isomorphic moduli space $\mathcal{A}_3(1,4,4)$. The result is based on the study the Hurwitz space $\mathcal{H}_{4,n}(Y)$ of 
 quadruple coverings of an elliptic curve $Y$ simply branched in $n\geq 2$ points. We prove the unirationality of its codimension one subvariety $\mathcal{H}^{0}_{4,A}(Y)$ which parametrizes quadruple coverings $\pi:X\to Y$ with Tschirnhausen modules isomorphic to $A^{-1}$, where  $A\in Pic^{n/2}Y$, and for which $\pi^*:J(Y)\to J(X)$ is injective. This  is an analog of the result of Arbarello and Cornalba that the Hurwitz space $\mathcal{H}_{4,n}(\mathbb{P}^1)$ is unirational.
\end{abstract}
\date{}
\author{Vassil Kanev}
\address{Department of Mathematics, University of Palermo\\
Via Archirafi n. 34, 90123 Palermo, Italy\\and\\
Institute of Mathematics, Bulgarian Academy of Sciences
}
\email{kanev@math.unipa.it}
\thanks{Research supported by the italian MIUR under the  program
"Geometria sulle variet\`{a} algebriche"}

\maketitle

\section*{Introduction}
In the present paper we study the Hurwitz space $\mathcal{H}_{4,n}(Y)$ 
which para\-metrizes simple quadruple coverings of an elliptic curve $Y$ 
branched in $n\geq 2$ points. There is a canonical smooth fibration 
$h:\mathcal{H}_{d,n}(Y)\to Pic^{n/2}Y$ (see (\ref{s2.8})). The fiber 
over $A$, which we denote by 
$\mathcal{H}_{d,A}(Y)$, parametrizes the coverings whose Tschirnhausen 
module has 
determinant isomorphic to $A^{-1}$. When studying coverings of an 
elliptic curve of non-prime degree it is natural to consider  
coverings which satisfy the condition that $\pi_*:H_1(X,\mathbb{Z})\to 
H_1(Y,\mathbb{Z})$ is surjective, or equivalently that $\pi^*:J(Y)\to J(X)$ is injective,  since coverings for which $Coker\, 
\pi_{*}\ne 0$ are reduced to coverings of smaller degree via  certain 
\'{e}tale coverings $\tilde{Y}\to Y$. We denote by 
$\mathcal{H}^{0}_{d,n}(Y)$ and $\mathcal{H}^{0}_{d,A}(Y)$ the 
corresponding  Hurwitz spaces.

One of the two main results of the paper states that $\mathcal{H}^0_{4,n}(Y)$ is connected and
$\mathcal{H}^{0}_{4,A}(Y)$ is connected and unirational 
(Theorem~\ref{s2.42}). We mention that $\mathcal{H}_{4,n}(Y)$ has three other connected components (see Remark~1 following Theorem~\ref{s2.42}).
In our previuos paper \cite{K1} we proved that 
$\mathcal{H}_{d,A}(Y)$ is unirational for $d\leq 3$. We notice the 
analogy  with the well-known result of 
Arbarello-Cornalba \cite{AC} which states that 
$\mathcal{H}_{d,n}(\mathbb{P}^{1})$ is unirational if $d\leq 5$. 
O. Schreyer gave in \cite{Sch} an alternative proof of the 
unirationality of $\mathcal{H}_{d,n}(\mathbb{P}^{1})$ for $d\leq 5$ by 
a method which was then developed by Casnati and Ekedahl in \cite{CE} 
(for $d=4$) and Casnati in \cite{Ca} (for $d=5$) and yields a 
description of Gorenstein coverings of degree 4 and 5 by means of 
a pair of vector bundles over the base, a connection analogous to the 
one found by Miranda for triple coverings \cite{Mi}. In our proof of 
the unirationality of $\mathcal{H}^{0}_{4,A}(Y)$ we use the result of 
\cite{CE}. Atiya's theory of vector bundles over elliptic curves \cite{At} is 
heavily used.

The way we prove the unirationality of $\mathcal{H}^{0}_{4,A}(Y)$ 
suggests a construction by which we are able to prove the 
unirationality of $\mathcal{A}_{3}(1,1,4)$, the moduli space of 
abelian threefolds with polarization of type $(1,1,4)$, and this is 
our second main result (Theorem~\ref{s3.58}). We follow the same 
pattern by which in \cite{K1} we proved the unirationality of 
$\mathcal{A}_{3}(1,1,d)$ for $d\leq 3$. We consider 
a smooth elliptic 
fibration $q:\mathcal{Y}\to B\subset \mathbb{P}^{1}$ and construct a 
pair of vector bundles $\mathcal{F}$ and $\mathcal{E}$ 
over $\mathcal{Y}$ by which we construct a family of quadruple 
coverings with a rational base $T$. One associates with such a family 
the Prym mapping $\varPhi:T\to \mathcal{A}_{3}(1,1,4)$ and we deduce 
the unirationality of $\mathcal{A}_{3}(1,1,4)$ by proving that 
$\varPhi$ is dominant. According to a theorem of Birkenhake-Lange 
\cite{BL3}\quad $\mathcal{A}_{3}(1,1,4)\cong \mathcal{A}_{3}(1,4,4)$. 
Thus our result establishes the unirationality of 
$\mathcal{A}_{3}(1,4,4)$ as well. We refer to the introduction of 
 \cite{GP}  the reader who may be interested about  other 
results on the unirationality of moduli spaces of abelian varieties. 
\begin{notation} 
A morphism (or 
holomorphic mapping) $\pi : X\to Y$ is called covering if it is finite, 
surjective and flat. 
Unless otherwise specified
 we make distinction between locally free sheaves and 
vector bundles and we denote differently their projectivizations. If $E$ is a 
locally free sheaf of $Y$ and if $\mathbb{E}$ is the corresponding vector 
bundle, i.e. $E\cong \mathcal{O}_Y(\mathbb{E})$, then $\mathbf{P}(E):=
\mathbf{Proj}(S(E)) \cong \mathbb{P}(\mathbb{E}^{\vee})$. 
Unless otherwise specified or clear from the context a curve is assumed to be irreducible. 
A covering of  projective curves $\pi : X\to Y$ of degree 
$d$ is called simple if $X$ and $Y$ are   
smooth and for each $y\in Y$ one has $d-1\leq \#\; \pi^{-1}(y) \leq d$. 
If $E$ is a locally free sheaf of a smooth curve $Y$ we denote the rank, the degree and the slope of $E$ by $r(E), d(E)$ and $\mu(E)=d(E)/r(E)$. If $Y$ is an elliptic curve 
 there is a unique, up to isomorphism, indecomposable locally free sheaf of rank $r$ and degree 0 which has nontrivial sections \cite{At}. We denote it by $F_r$.
Unless 
otherwise specified we assume the base field $k=\mathbb{C}$.
\end{notation}

\section{Preliminaries}\label{s1}
Let $\pi :X\to Y$ be a finite covering of smooth, projective curves of 
degree $d\geq 2$, suppose $g(Y)\geq 1$. Let $P= 
Ker(Nm_{\pi}:J(X)\to J(Y))^{0} $ be the Prym variety of the covering. 
Let $\Theta$ be the canonical polarization of $J(X)$ and let 
$\Theta_{P}$ be its restriction on $P$. The following statement is proved in 
\cite{K1} Lemma~1.1.
\begin{pro}\label{s1.1}
The following three conditions are equivalent: $\pi_{*}: \linebreak
H_1(X,\mathbb{Z})\to 
H_1(Y,\mathbb{Z})$ is surjective; $\pi^{*}:J(Y)\to J(X)$ is injective; 
$Ker\, Nm_{\pi }$ is connected. Suppose these conditions hold and let 
$P=Ker\, Nm_{\pi }$. Then the polarization $\Theta_{P}$ is of type 
$(1,\ldots,1,d,\ldots,d)$ where the $d$'s are repeated $g(Y)$ times. 
\end{pro}

\begin{block}\label{s1.2}
Let $\pi :X\to Y$ be a covering as above, suppose $Y$ is an elliptic 
curve and $g(X)\geq 3$. Let 
$d_{2}=|H_1(Y,\mathbb{Z}):\pi_{*}H_1(X,\mathbb{Z})|$. Then $\pi$ may 
be decomposed as 
$X\overset{\pi_{1}}{\to}\tilde{Y}\overset{\pi_{2}}{\to}Y$ where 
$\pi_{2}$ is an isogeny of degree $d_{2}$ and $\deg \pi_{1} = d_{1} =
\frac{d}{d_{2}}$. According to the preceeding proposition the type of 
the polarization $\Theta_{P}$ is $(1,\ldots ,1,d_{1})$. In particular 
if $d=4$ one obtains polarization of type $(1,\ldots ,1,4)$ when 
$\pi_*:H_1(X,\mathbb{Z})\to H_1(Y,\mathbb{Z})$ is surjective and 
$(1,\ldots ,1,2)$ when it is not. Counting parameters we see that one 
might obtain a general abelian variety $A$ with polarization of type 
$(1,\ldots ,1,4)$ as Prym variety of a covering $\pi :X\to Y$ only if 
$\dim A=2$ or $\dim A = 3$. 
\end{block}

\begin{block}\label{s2.37}
We will need some facts about vector bundles over elliptic curves, a 
theory due to Atiyah \cite{At}. 
In this section we will make the customary identification between 
vector bundles and locally free sheaves.
For generalities on vector bundles 
over smooth, projective curves we refer the reader to the survey 
atricles \cite{Oe} and \cite{Br}. The following facts about vector 
bundles over an elliptic curve are known  (see e.g. \cite{Br} p.87). Every indecomposable vector 
bundle is semistable. A vector bundle is 
semistable if and only if it is a direct sum of indecomposable vector 
bundles of the same slope. A vector bundle of degree $d$ and rank $r$ 
is stable if and only if it is indecomposable and $(d,r)=1$. 
If $G$ is a semistable vector bundle over an 
elliptic curve, then $h^1(G)=0$ if $\mu(G)=\frac{d(G)}{r(G)}>0$ and 
$G$ is generated by its global sections if $\mu(G)>1$ (see e.g. \cite{Oe} p.39).
Let $r\geq 1$ and $d\in \mathbb{Z}$. The isomorphism classes of indecomposable vector bundles of rank $r$ and degree $d$ over an elliptic curve $Y$ depend on one parameter and are parametrized by $Pic^0(Y)\cong Y$. In fact if one fixes one such $E$, then all others are obtained as $E\otimes L$ for some $L\in Pic^0(Y) = J$. Furtermore $E\otimes L \cong E\otimes L'$ if and only if $L'\otimes L^{-1}$ is a point of order $r'=r/h$ of $Pic^0(Y)$ where $h=(r,d)$ (cf. \cite{At} Theorem~10). So, the indecomposable vector bundles of rank $r$ and degee $d$ are parametrized by $J/J_{r'}\cong J$. 

One may construct a Poincar\'{e} vector bundle as follows. 
Let us fix a point $y_{0}\in Y$.  Let $\mathcal{L}$ be the Poincar\'{e} line bundle on $Y\times J$ normalized by $\mathcal{L}|_{y_{0}\times J}\cong 
\mathcal{O}_{J}$. First suppose $(r,d)=1$. Then the general theory of 
stable vector bundles yields a Poincar\'{e} vector bundle 
$\mathcal{E}(r,d)$ on $Y\times U(r,d)$ where $U(r,d)$ is the fine 
moduli space of stable vector bundles (\cite{Ne} Ch.5~\S5). By \cite{At} p.~434 the 
morphism $U(r,d)\to J$ given by $u\mapsto \det 
\mathcal{E}(r,d)_{u}\otimes \mathcal{O}_{Y}(-dy_{0})$ is an 
isomorphism. So we may replace $U(r,d)$ by $J$. If $r=r'h, d=d'h$ 
where $(r',d')=1$ we let $\mathcal{E}(r,d) = \mathcal{E}(r',d')\otimes 
\pi_{1}^{*}F_{h}$. According to \cite{At} Lemma~24 and Theorem~10 
this family has the property that each indecomposable locally free 
sheaf $E$ of rank $r$ and degree $d$ is isomorphic to 
$\mathcal{E}(r,d)_{u}$ for some $u\in J$ and if $u\ne u'$, then 
$\mathcal{E}(r,d)_{u}\ncong \mathcal{E}(r,d)_{u'}$. Furthermore the 
invertible sheaf $\det \mathcal{E}(r,d)$ determines a morphism $J\to 
J$ given by $u\mapsto \det \mathcal{E}(r,d)_{u}\otimes 
\mathcal{O}_{Y}(-dy_{0})$ and this morphism coincides with $h\cdot 
id_{J}:J\to J$. 

The folllowing statement is proved in \cite{Te} Lemma~2.3.
\end{block}

\begin{pro}\label{s1.3}
Let $E$ and $G$ be semistable vector bundles over an elliptic curve. Then 
$E\otimes G$ is semistable of slope $\mu(E\otimes G)=\mu(E)+\mu(G)$.
\end{pro}

\begin{cor}\label{s1.5}
Let $E_{1}, E_{2}, E_{3}, E_{4}, \ldots$ be semistable vector bundles 
over an elliptic curve. Then
\begin{equation}\label{es1.5}
S^{k_{1}}E_{1}\otimes \wedge^{k_{2}}E_{2}\otimes 
S^{k_{3}}E_{3}^{*}\otimes \wedge^{k_{4}}E_{4}^{*}\otimes \cdots
\end{equation}
is semistable of slope $\mu = 
k_{1}\mu(E_{1})+k_{2}\mu(E_{2})-k_{3}\mu(E_{3})-k_{4}\mu(E_{4})+\cdots$
\end{cor}
\begin{proof}
By Proposition~\ref{s1.3} the tensor product 
$E_{1}^{\otimes k_{1}}\otimes E_{2}^{\otimes k_{2}}\otimes 
(E_{3}^*)^{\otimes k_{3}}\otimes (E_{4}^*)^{\otimes k_{4}}\otimes \cdots
$ is semistable, so it is a direct sum of indecomposable vector 
bundles of slope $\mu$. The vector bundle \eqref{es1.5} is its direct 
summand, so the corollary follows from an analog of the Krull-Schmidt theorem 
due to Atiyah \cite{At1}.
\end{proof}

\bigskip
\noindent
If $E$ is a vector bundle of rank 2 and even degree over an elliptic 
curve the symmetric powers $S^{n}E$ were calculated by Atiyah: if 
$E\cong F_{2}\otimes L$, then $S^{n}(F_{2}\otimes L)\cong F_{n}\otimes 
L^{n}$ (cf. \cite{At} p.438). In the case of rank 2 vector bundles of odd 
degree the following statement holds.

\begin{pro}\label{s1.5a}
Let $E$ be  a vector bundle over an elliptic curve of rank 2 and odd 
degree. Then for every $k\geq 1$ one has
\renewcommand{\theenumi}{\roman{enumi}}
\begin{enumerate}
\item $S^{2k-1}E\cong (\wedge^{2}E)^{k-1}\otimes (E^{\oplus k})$,
\item $S^{2k}E\cong (\wedge^{2}E)^{k}\otimes 
\left[I^{\oplus(2k+1-3
\lceil \frac{k}{2}\rceil)}\oplus(\oplus_{i=1}^{3}L_{i})^{\oplus 
\lceil \frac{k}{2}\rceil}\right]$
\end{enumerate}
where $L_{i}^{2}\cong I, L_{i}\ncong I$.
\end{pro}
\begin{proof}
We prove the formula for $S^{n}E$ by induction on $n$. Let $n=2$. Let 
$A$ be a fixed line bundle of degree 1 as in \cite{At}. One has 
$E\cong E_A(2,1)\otimes M$ for some line bundle $M$. According to 
Lemma 22 (ibid.) $E\otimes E^{*}\cong \oplus_{i=0}^{3}L_{i}$, where 
$L_{0}=I$. Since $E\cong E^{*}\otimes\wedge^{2}E$ one obtains $E\otimes 
E\cong \wedge^{2}E\otimes (\oplus_{i=1}^{3}L_{i})$. On the other hand 
$E\otimes E\cong S^{2}E\oplus \wedge^{2}E$. So $S^{2}E\cong 
\wedge^{2}E\otimes \sum_{i=1}^{3}L_{i}$ by  
\cite{At1}. Let us prove the formula for $2k+1, 2k+2$ assuming it 
holds for $n\leq 2k$. By the Clebsch-Gordan formula $S^{2k}E\otimes E\cong 
S^{2k+1}E\oplus (\wedge^{2}E\otimes S^{2k-1}E)$ which is isomorphic to 
$S^{2k+1}E\oplus ((\wedge^{2}E)^{k}\otimes E^{\oplus k})$ by the induction 
hypothesis applied to $n=2k-1$. One has $E\otimes L_{i}\cong E$ (cf. \cite{At} 
p.434), so by the induction hypothesis applied to $n=2k$ one has 
$S^{2k}E\otimes E \cong (\wedge^{2}E)^{k}\otimes E^{\oplus (2k+1)}$. 
Using  \cite{At1} one obtains
\(
S^{2k+1}E\ \cong \ (\wedge^{2}E)^{k}\otimes E^{\oplus (k+1)}.
\)
Let $a_{k}=2k+1-3\lceil \frac{k}{2}\rceil, b_{k}=\lceil 
\frac{k}{2}\rceil$. One has 
\begin{multline*}
S^{2k+1}E\otimes E \cong S^{2k+2}E\oplus (\wedge^{2}E\otimes 
S^{2k}E)\\
\cong S^{2k+2}E \oplus (\wedge^{2}E)^{k+1}\left[I^{a_{k}}\oplus 
(\sum_{i=1}^{3}L_{i})^{b_{k}}\right].
\end{multline*}
On the other hand
\begin{multline*}
S^{2k+1}E\otimes E \cong ((\wedge^{2}E)^{k}\otimes E^{\oplus 
(k+1)})\otimes E \\
\cong (\wedge^{2}E)^{k}\otimes (E\otimes E)^{\oplus 
(k+1)} \cong 
(\wedge^{2}E)^{k+1}\otimes \left[I^{k+1}\oplus(\sum_{i=1}^{3}L_{i})^{k+1}\right].
\end{multline*}
By \cite{At1} we may cancel and obtain 
\[
S^{2k+2}E\ \cong \ (\wedge^{2}E)^{k+1}\otimes \left[I^{k+1-a_{k}}\oplus 
(\sum_{i=1}^{3}L_{i})^{k+1-b_{k}}\right].
\]
That $k+1-a_{k}=a_{k+1}$ and $k+1-b_{k}=b_{k+1}$ is clear. 
\end{proof}

\section{Hurwitz spaces of quadruple coverings of elliptic \\ curves}\label{s2}
\begin{block}\label{s2.6}
We study the Hurwitz spaces of quadruple coverings by means of a 
result of Casnati and Ekedahl \cite{CE} which describes such coverings 
in terms of a pair of vector bundles on the base. We recall their 
construction in the special case we need. Let $\pi :X\to Y$ be a 
covering of smooth projective curves of degree $d$. The Tschirnhausen 
module of the covering is the quotient sheaf $E^{\vee}$ defined by the 
exact sequence
\[
0\lto 
\mathcal{O}_{Y}\overset{\pi^{\#}}{\lto}\pi_{*}\mathcal{O}_{X}\lto 
E^{\vee}\lto 0.
\]
One has $E^{\vee}\cong Ker(Tr_{\pi}:\pi_{*}\mathcal{O}_{X}\to 
\mathcal{O}_{Y})$ and this is a locally free sheaf of rank $d-1$. There 
is a canonical embedding $i:X\to \mathbf{P}(E)$ defined by the relative 
dualizing sheaf $\omega_{X/Y}\cong 
\omega_{X}\otimes (\pi^{*}\omega_{Y})^{-1}$ and satisfying 
$i^{*}\mathcal{O}_{\mathbf{P}(E)}(1)\cong \omega_{X/Y}$.
When $d=4$ every fiber $X_{y}$ is an intersection of two conics 
in $\mathbf{P}(E)_{y}$. This globalizes as follows. There is a locally 
free sheaf $F$ of rank 2 on $Y$ such that if $\rho:\mathbf{P}(E)\to Y$ 
is the canonical fibration and $N=\rho^{*}F$, then the resolution of 
$\mathcal{O}_{X}$ is given by 
\[
0\to \rho^{*}\det E(-4) \to N(-2) \overset{\delta}{\to} 
\mathcal{O}_{\mathbf{P}(E)}\to \mathcal{O}_{X}\to 0.
\]
By the results of Casnati and Ekedahl one has $\det F\cong \det E$ and 
given such a pair of locally free sheaves on $Y$ a quadruple covering 
is determined uniquely by the homomorphism $\delta \in 
H^0(\mathbf{P}(E),\Check{N}(2))$. Let 
$\phi:H^0(Y,\Check{F}\otimes S^{2}E)\to 
H^0(\mathbf{P}(E),\Check{N}(2))$ be the canonical isomorphism and let 
$\delta = \phi(\eta)$. Then $\eta$ satisfies the following 
property. For every $y\in Y$ the value $\eta(y):F(y)\to S^{2}E(y)$ 
determines a pencil of conics in $\mathbb{P}^{2}=\mathbf{P}(E)_{y}$ 
whose base locus is of dimension 0. According to \cite{CE} Definition~4.2 the sections $\eta$ with this property are called \emph{of right 
codimension for every $y\in Y$}. The described construction is valid 
for every integral scheme $Y$ and yields a bijective correspondence 
between the following data (\cite{CE} Theorem~4.4.)
\renewcommand{\theenumi}{\Alph{enumi}}
{\it 
\begin{enumerate}
\item Finite, flat Gorenstein coverings $\pi :X\to Y$ of degree 4.
\item Locally free sheaves $F$ of rank 2 and $E$ of rank 3  on $Y$ 
such that $\det F \cong \det E$, and a section $\eta \in 
H^0(Y,\Check{F}\otimes S^2E)$ which has right codimension at every 
$y\in Y$. The section is determined uniquely up to multiplication by 
scalars.
\end{enumerate}
}
Given $F$ and $E$ with $\det F \cong \det E$ one can prove that the 
subset of $\mathbb{P}H^0(Y,\Check{F}\otimes S^2E)$ which parametrizes 
the elements $\langle \eta \rangle$  which have right codimension at 
every $y\in Y$ is Zariski open (cf. Lemma~\ref{s2.34a}). In case it 
is not empty, by \cite{CE} Theorem~4.4, the group $PGL(F)\times PGL(E)$ 
acts faithfully on this Zariski open set. The orbits of this 
action correspond bijectively to 
the equivalence classes of Gorenstein coverings of degree 4 of $Y$ 
whose canonically associated pair of locally free sheaves is 
isomorphic to the pair $(E,F)$.
\end{block}

\begin{block}\label{s2.8}
Let $Y$ be a smooth, projective curve. Let 
$\mathcal{H}_{d,n}(Y)$ be the Hurwitz space which para\-metrizes the 
simple coverings of $Y$ of degree $d$ branched in $n$ points. We 
denote by $\mathcal{H}^{0}_{d,n}(Y)$ the subset of the Hurwitz space
$\mathcal{H}_{d,n}(Y)$ whose points correspond to coverings $\pi :X\to 
Y$  with the property that $\pi_*:H_1(X,\mathbb{Z})\to 
H_1(Y,\mathbb{Z})$ is surjective. This property is preserved under 
deformation, so $\mathcal{H}^{0}_{d,n}(Y)$ is a union of connected 
components of $\mathcal{H}_{d,n}(Y)$. Given $A\in Pic^{n/2}(Y)$ we 
denote by $\mathcal{H}_{d,A}(Y)$ the closed reduced subscheme whose 
points correspond to coverings having a Tschirnhausen module with 
determinant isomorphic to $A^{-1}$ (cf. \cite{K1} (2.4) and Lemma~2.5). We denote by $\mathcal{H}^{0}_{d,A}(Y)$ the intersection 
$\mathcal{H}^{0}_{d,A}(Y)\cap \mathcal{H}_{d,A}(Y)$. In this paper we 
study  $\mathcal{H}^{0}_{4,n}(Y)$ and 
$\mathcal{H}^{0}_{4,A}(Y)$ where $Y$ is elliptic curve and $n\geq 2$. 
Notice that both Hurwitz spaces are nonempty. Indeed, if 
$\pi_*:H_1(X,\mathbb{Z})\to H_1(Y,\mathbb{Z})$ is not surjective, then 
$\pi$ may be decomposed as $X\to \tilde{Y}\to Y$, where $\tilde{Y}\to 
Y$ is an unramified double covering. Hence the monodromy group of $\pi 
:X\to Y$ is different from $S_{4}$. Simple branched coverings with 
 monodromy group $S_{4}$ are easily constructed (see e.g. \cite{K1} 
Lemma~2.1), so $\mathcal{H}^{0}_{4,n}(Y)\ne \emptyset$ for every pair 
$n\geq 2$. Using Lemma~2.5 of \cite{K1} one concludes that 
$\mathcal{H}^{0}_{4,A}(Y)\ne \emptyset$ for every $A\in Pic^{n/2}Y$.
\end{block}

\noindent
The unirationality results proved in the paper are based on 
the following theorem.

\begin{thm}\label{s2.9}
Let $Y$ be an elliptic curve. Let $n=2e\geq 2$. There is a Zariski 
open, dense subset $U\subset \mathcal{H}^{0}_{4,n}(Y)$ such that for 
every $[X\to Y]\in U$ one has for the associated pair of locally free 
sheaves $(E,F)$:
\begin{enumerate}
\item If $e \not \equiv 0 (mod\, 3)$ then $E$ is indecomposable of 
degree $e$.
\item If $e \equiv 0 (mod\, 3)$, then $E$  is isomorphic to 
$L_{1}\oplus L_{2}\oplus L_{3}$ where $L_{i}\in Pic^{e/3}Y$ and 
$L_{i}\ncong L_{j}$ for $i\ne j$.
\item If $e \equiv 1(mod\, 2)$ then $F$ is indecomposable of degree $e$. 
\item If $e \equiv 0(mod\, 2)$, then $F\cong M_{1}\oplus M_{2}$ where 
$M_{i}\in Pic^{e/2}Y$ and $M_{1}\ncong M_{2}$.
\end{enumerate}
\end{thm}
The proof of the theorem is rather long, it occupies (\ref{s2.11a}) -- 
(\ref{s2.38}), and we indicate first the main steps. We use the 
result of Casnati and Ekedahl \cite{CE} by which quadruple coverings 
of an elliptic curve are described in terms of a pair of vector 
bundles $(E,F)$ of ranks 3 and 2 respectively. Our fisrt task is to 
bound from 
above the number of parameters on which the coverings depend when the 
pair varies within a family of a given type. This quite long 
calculation is done in 7 steps according to the polygons of the 
Harder-Narasimhan filtrations of the vector bundles $E$ and $F$. The outcome 
is that only the types specified in the theorem may give sufficient 
number of moduli. This work is done in (\ref{s2.11a}) -- 
(\ref{s2.32}). The remaining part of the proof is devoted to show that 
for each of the occuring types of $(E,F)$ one may construct a family of 
quadruple coverings $\mathcal{X}\to Y\times T$ such that when applying 
the universal property of the Hurwitz space one obtains a morphism 
$f:T\to \mathcal{H}_{4,n}(Y)$ such that the dimension of the image of 
$T$ equals the number of parameters calculated  in the first part of 
the proof. Here the technical result of
 Lemma~\ref{s2.34a}  is used.  This is also applied  
later in connection with the moduli space $\mathcal{A}_{3}(1,1,4)$. This lemma together with the calculation of the number of parameters permits to conclude in (\ref{s2.38}) that every sufficiently general quadruple covering of $\mathcal{H}^0_{4,n}(Y)$ has a pair $(E,F)$ of  the type specified in the theorem.

\begin{lem}\label{s2.11a}
Let $E$ and $F$ be locally free sheaves over an elliptic curve of 
ranks $r(E)=3, r(F)=2$. Let $\deg E = \deg F = e$. Then 
\[
\chi(\Check{F}\otimes S^2E) = h^0(\Check{F}\otimes S^2E) - 
h^1(\Check{F}\otimes S^2E) = 2e
\]
\end{lem}
\begin{proof}
By Riemann-Roch $\chi(\Check{F}\otimes S^2E)=\deg(\Check{F}\otimes 
S^2E)$. One has $c_{1}(\Check{F}\otimes S^2E)=8c_{1}(E)-6c_{1}(F)$. 
Thus $\deg(\Check{F}\otimes S^2E) = 2e$ .
\end{proof}

\begin{block}\label{s2.11}
For the various types of $(E,F)$ considered below the following 
numbers are constant: $h^0(End\, E),\: h^0(End\, F),\: 
h^0(\Check{F}\otimes S^2E)$. Specifying the type of a pair $(E,F)$ means generally specifying the decomposition types of $E$ and $F$, the possibility about isomorphy or non isomorphy of various summands that appear, and when $h^0(\Check{F}\otimes S^2E)$ may jump, imposing additional conditions on $F$ and $E$ in order that $h^0(\Check{F}\otimes S^2E)$ stays fixed. By (\ref{s2.6}) given a pair $(E,F)$, if 
there are Gorenstein coverings of degree 4 whose associated pair is 
isomorphic to the given one, then such coverings depend on the 
following number of parameters:
\begin{multline*}
 \dim 
\mathbb{P}H^0(\Check{F}\otimes S^2E) - \dim PGL(E) -\dim PGL(F)\\ 
 =\ 
h^0(\Check{F}\otimes S^2E) - h^0(End\, E) - h^0(End\, F) +1
\end{multline*}
Varying $(E,F)$ in a family of a given type one obtains
\begin{multline}\label{es2.11}
\# \, \text{moduli} \, [X\to Y]\ =\\
 \# \, \text{moduli} \, (E,F) +
h^0(\Check{F}\otimes S^2E) - h^0(End\, E) - h^0(End\, F) + 1.
\end{multline}
Let us consider the Hurwitz space $\mathcal{H}_{4,n}(Y)$ parametrizing 
equivalence classes of simple coverings of an elliptic curve $\pi 
:X\to Y$ branched in $n$ points. Here the Tschirnhausen modules have 
degree $-e$ where $n=2e$ (\cite{K1} Lemma~2.3). Hence by 
Lemma~\ref{s2.11a} for the pairs $(E,F)$ associated to the points of 
$\mathcal{H}_{4,n}(Y)$ one has $h^0(\Check{F}\otimes S^2E) = n + 
h^1(\Check{F}\otimes S^2E)$. Taking into account that $h^0(End\, E) = 
h^1(End\, E), \: h^0(End\, F) = h^1(End\, F)$ we obtain the following 
criterion:

\noindent
\emph{Pairs of locally free sheaves $(E,F)$ of a given type yield 
quadruple coverings $\pi :X\to Y$  with insufficient number of 
moduli, i.e. $\# \, \text{moduli} \, [X\to Y] < n = \dim 
\mathcal{H}_{4,n}(Y)$ if}
\begin{equation}\label{es2.12}
\begin{split}
\# \, \text{moduli} \, (E,F) 
< &h^1(End\, E) + h^1(End\, F) - 1 - h^1(\Check{F}\otimes S^2E) = \\
  &h^0(End\, E) + h^0(End\, F) - 1 - h^1(\Check{F}\otimes S^2E)
\end{split}
\end{equation}
\end{block}

\begin{block}\label{s2.12}
 {\bf Case 1.} \emph{$E$ and $F$ are semistable locally free sheaves}. 
According to Corollary~\ref{s1.5}\quad $\Check{F}\otimes S^2E$ is 
semistable. Its slope $\mu(\Check{F}\otimes S^2E) = -\mu(F) + 
2\mu(E)=-\frac{e}{2}+2\cdot\frac{e}{3}=\frac{e}{6}>0$. Hence 
$h^1(\Check{F}\otimes S^2E)=0$. Let us decompose in a direct sum of 
indecomposable locally free sheaves: $E=E_{1}\oplus \cdots \oplus 
E_{\ell},\: F=G_{1}\oplus \cdots \oplus G_{m}$. One has 
$\mu(E_{i})=\mu(E),\: \mu(G_{j})=\mu(F)$ for every $i,j$. One obtains 
various possible types of the pairs $(E,F)$ by: a) fixing the ranks of 
$E_{i}$ and $G_{j}$; b) requiring that certain direct summands of $E$ 
are isomorphic to each other and similarly for $F$. We have $\# \, 
\text{moduli} \, (E,F) \leq \ell +m -1$. Here subtracting 1 is for 
$\det E \cong \det F$. One has a strict inequality if one considers a 
type where certain direct summands are isomorphic to each other. So, 
in this case the inequality \eqref{es2.12} holds since $\ell \leq 
h^1(End\, E)$ and $m\leq h^1(End\, F)$. Suppose one considers a type 
of $(E,F)$ where one of the direct 
summands $E_{i}$ is indecomposable of rank $r'h$ and degree $d'h$ 
where $h>1$ and $(r',d')=1$. Then $h^1(End\, E_{i})=h$ (cf. the proof 
of \cite{At} Lemma~23). Therefore $\ell < h^1(End\, E)$ and 
\eqref{es2.12} holds. A similar argument may be applied to $F$. 
If every direct summand $E_{i},\: G_{j}$ has degree and rank prime to 
each other, then $r(E_{i})=r(E_{j})$ and $r(G_{i})=r(G_{j})$ for $i\ne 
j$. Indeed $\mu(E_{i})=\mu(E_{j})$ implies 
$d(E_{i})r(E_{j})=d(E_{j})r(E_{i})$, so $r(E_{i})|r(E_{j})$ and 
$d(E_{i})|d(E_{j})$. This implies $r(E_{i})=r(E_{j}),\: 
d(E_{i})=d(E_{j})$. We conclude that the only possible types of $(E,F)$ 
which might give $\# \, \text{moduli} \, [X\to Y] = n$ are the types 
$E\cong \oplus_{i=1}^{\ell}E_{i},\; F\cong \oplus_{j=1}^{m}G_{j}$ 
where each $E_{i}$ or $G_{j}$ is indecomposable with rank prime to its 
degree, $r(E_{i})=r(E_{j}),\; r(G_{i})=r(G_{j})$ and furthermore 
$E_{i}\ncong E_{j}$ and $G_{i}\ncong G_{j}$ for $i\ne j$. When $n=2e$ 
is fixed there is only one such type and this is the type of 
Theorem~\ref{s2.9}. All other possible types with semistable $E$ and 
$F$ yield $\# \, \text{moduli} \, [X\to Y]<n$ according to the 
criterion of (\ref{s2.11}).
\end{block}

\begin{block}\label{s2.14}
We need an explicit form of $X\subset \mathbf{P}(E)$ in order to 
exclude some types of $(E,F)$. Let $U\subset Y$ be a Zariski open 
subset such that  $E|_U$ and $F|_U$ are trivial. Let 
$E|_{U}=\mathcal{O}_{U}e_{1}\oplus \mathcal{O}_{U}e_{2}\oplus 
\mathcal{O}_{U}e_{3},\quad 
F|_{U}=\mathcal{O}_{U}f_{1}\oplus \mathcal{O}_{U}f_{2}$. One has 
$S^{2}E|_{U}=\sum_{i\leq j}\mathcal{O}_{u}e_{i}e_{j}$. If $\eta \in 
H^{0}(\Check{F}\otimes S^2E) = Hom_{Y}(F,S^{2}E)$, then
\begin{equation*}
\begin{split}
\eta(\alpha _{1}f_{1}+\alpha _{2}f_{2}) 
&= \alpha _{1}\eta(f_{1}) + \alpha _{2}\eta(f_{2})\\
&= \alpha _{1}\sum_{i\leq j}a_{ij,1}e_{i}e_{j} + \alpha 
_{2}\sum_{i\leq j}a_{ij,2}e_{i}e_{j}
\end{split}
\end{equation*}
If $f^{1},f^{2}$ are the sections of $\Check{F}|_{U}$ dual to 
$f_{1},f_{2}$, i.e. $\langle f_{i},f^{j}\rangle=\delta_{ij}$, then 
$\eta \in H^{0}(\Check{F}\otimes S^2E)$ is locally given by
\begin{equation}\label{es2.14}
\eta|_{U}\ =\ \sum_{i\leq j}a_{ij,1}f^{1}\otimes e_{i}e_{j}
+ \sum_{i\leq j}a_{ij,2}f^{2}\otimes e_{i}e_{j}.
\end{equation}
For every $y\in U$ the fiber $\mathbf{P}(E)_{y}\cong \mathbb{P}^{2}$ 
has coordinates $e_{1}(y),e_{2}(y),e_{3}(y)$. Every pair 
$\xi_{1},\xi_{2}\in \mathbb{C}$ with $(\xi_{1},\xi_{2})\ne (0,0)$ 
defines a conic with equation
\begin{equation}\label{es2.14a}
\sum_{i\leq j}(\xi_{1}a_{ij,1}(y)+\xi_{2}a_{ij,2}(y))e_{i}(y)e_{j}(y)\ 
=\ 0.
\end{equation}
Varying $(\xi_{1},\xi_{2})$ one obtains a pencil of conics whoce base 
locus is $X_{y}\subset \mathbf{P}(E)_{y}$.
\end{block}

\begin{block}\label{s2.15}
{\bf Case 2.} \emph{$F$ is semistable, $E\cong E_{1}\oplus E_{2}$, where $E_{i}$ 
are semistable of slopes $\mu(E_{1})=\mu_{1}>\mu(E_{2})=\mu_{2}$}. 
This case is subdivided further according to $r(E_{1})$, but we treat 
the two subcases simultaneously whenever possible. We have 
$S^{2}E\cong S^{2}E_{1}\oplus (E_{1}\otimes E_{2})\oplus S^{2}E_{2}$. 
Here $S^{2}E_{2}=E_{2}^{2}$ if $r(E_{2})=1$ and 
$S^{2}E_{1}=E_{1}^{2}$ if $r(E_{1})=1$.

{\bf Subcase 2A.}\quad $\mu(\Check{F}\otimes 
S^2E_{2})=-\frac{e}{2}+2\mu_{2}>0$, or $\mu(\Check{F}\otimes 
S^2E_{2})=0$ and $h^1(\Check{F}\otimes S^2E_{2})=
h^0(\Check{F}\otimes S^2E_{2})
=0$. Here
\[
\mu(\Check{F}\otimes 
S^2E_{1})>\mu(\Check{F}\otimes E_{1}\otimes 
E_{2}) > \mu(\Check{F}\otimes S^2E_{2}) \geq 0
\]
so we obtain $h^1(\Check{F}\otimes S^2E)=0$. It is clear that 
\[
\# \, \text{moduli} \, (E,F) \leq h^1(End\, E_{1}) + h^1(End\, 
E_{2}) + h^1(End\, F) - 1
\]
We have $h^1(End\, E) = h^1(End\, E_{1})+h^1(End\, E_{2}) + 
h^1(E_{1}^{*}\otimes E_{2})$ and $h^1(E_{1}^{*}\otimes E_{2}) = h^0(
E_{1}\otimes E_{2}^{*})= 2(\mu_{1}-\mu_{2}) > 0$. We conclude 
inequality \eqref{es2.12} holds whatever type of $(E,F)$ satisfying 
the conditions of this subcase is considered.

{\bf Subcase 2B.}\quad $\mu(\Check{F}\otimes S^2E_{2})=-\frac{e}{2}+2\mu_{2}<0$. 
Here $H^{0}(\Check{F}\otimes S^2E_{2})=0$. If $r(E_{2})=1$ this means 
that if we choose in (\ref{s2.14}) the frame of $E|_{U}$ so that 
$e_{1},e_{2}$ generate $E_{1}|_{U}$ and $e_{3}$ generates 
$E_{2}|_{U}$, then $a_{33,1}=a_{33,2}=0$. This implies that in each 
fiber $\mathbf{P}(E)_{y},\; y\in U$ the point $e_{1}(y)=e_{2}(y)=0$ 
belongs to the pencil of conics defined by \eqref{es2.14a}. This means 
that the section of $\mathbf{P}(E)\to Y$ defined by $E\to E_{1}\to 
0$ is a component of the Gorenstein covering $\pi :X\to Y$ determined 
by any $\eta\in H^0(\Check{F}\otimes S^2E)$. We may thus exclude from 
consideration such type, since no irreducible, smooth, quadruple 
cover $X$ may be obtained from such a pair $(E,F)$. If $r(E_{1})=1, 
r(E_{2})=2$ choosing $e_{i}$ in (\ref{s2.14}) so that $e_{1}$ 
generates $E_{1}|_{U}$ and $e_{2},e_{3}$ generate $E_{2}|U$ we see 
that for $\forall y\in U$ all conics of \eqref{es2.14a} must contain 
the line $\{e_{1}(y)=0\}$. So no pair $(E,F)$ with this property could 
be associated with a Gorenstein covering $X\to Y$.

{\bf Subcase 2C.}\quad  $\mu(\Check{F}\otimes S^2E_{2})=0,\; h^0(\Check{F}\otimes 
S^2E)>0$. Here $h^1(\Check{F}\otimes S^2E)=h^0(\Check{F}\otimes 
S^2E_{2})$. We have $-\mu(F)+2\mu(E_{2})=0$, so $\mu(F)$ is an integer 
and $e=d(F)=2\mu(F)$ is even. One has two possibilities. Either $F$ is 
indecomposable of even degree or $F\cong L_{1}\oplus L_{2}$ where 
$\deg L_{i}=\frac{e}{2}$. We recall $F_r$ denotes the unique (up to isomorphism) locally free sheaf on $Y$ of rank $r$ and degree $0$ with $h^0(F_r)\geq 1$.

{\bf Subcase 2C$'$.}\quad Assume $F$ is indecomposable, $F\cong 
F_{2}\otimes L$ with $\deg L = e/2$. In this case the equality $\det 
F\cong \det E$ determines $F$ by $E$ up to tensoring by an element of 
$(Pic^{0}Y)_{2}$. Thus $\#\, \text{moduli}\, (E,F) = \#\, 
\text{moduli}\, (E)$. One has $h^1(End\, F) = h^0(End 
F_{2})=h^0(F_{1}\oplus F_{2}) = 2$ (cf. \cite{At} p.437) and furthermore
\begin{equation}\label{es2.16}
\#\, \text{moduli}\, (E) \leq h^1(End\, E_{1}) + h^1(End\, E_{2})
\end{equation}
So one would have the inequality \eqref{es2.12} for any type where it 
holds
\begin{equation}\label{es2.17}
1 + h^0(E_{1}\otimes E_{2}^{*}) - h^0(\Check{F}\otimes S^2E) > 0
\end{equation}
We have $\mu(E_1\otimes E_2^*) = \mu_{1}-\mu_{2}>0$. So $h^0(E_1\otimes 
E_2^*)=2(\mu_{1}-\mu_{2})\geq 1$. If $E_{2}$ is of rank 2 and 
decomposable then $\mu_{2}$ is an integer, so $h^0(E_1\otimes 
E_2^*)\geq 2$. If $E_{2}$ is of rank 1, then $\Check{F}\otimes 
S^2E_{2}=\Check{F}\otimes E^{2}_{2}$ is indecomposable of degree 0, so 
$h^0(\Check{F}\otimes E_{2}^{2})\leq 1$. This shows \eqref{es2.17} 
holds. If $E_{2}$ is of rank 2 one has the following cases.

{\bf Subcase 2C$'$(i).}\emph{\quad $E_{2}$ is indecomposable of odd degree.} 
Then according to Proposition~\ref{s1.5a} one has 
\(
S^{2}E_{2}\cong \wedge^{2}E_{2}\otimes (\eta_{1}\oplus\eta_{2}\oplus 
\eta_{3})\) where \(\eta_{i}^{2}\cong \mathcal{O}_{Y},\; 
\eta_{i}\ncong \mathcal{O}_{Y}.
\)
Let $F\cong F_{2}\otimes L$. Then $\Check{F}\otimes S^2E_{2}$ is a 
direct sum of $F_{2}\otimes L^{-1}\otimes \wedge^{2}E_{2}\otimes 
\eta_{i},\; i=1,2,3$. According to \cite{At} Theorem~5 only one of 
these locally free sheaves might be isomorphic to $F_{2}$. Thus 
$h^0(\Check{F}\otimes S^2E_{2})\leq 1$ and therefore \eqref{es2.17} 
holds.

{\bf Subcase 2C$'$(ii).}\quad \emph{$E_{2}$ is indecomposable of even 
degree.} Then $E_{2}\cong F_{2}\otimes M$ for some invertible sheaf 
$M$. Here $S^{2}E_{2}\cong F_{3}\otimes M^{2}$ and \linebreak
$\Check{F}\otimes 
S^2E_{2}\cong F_{2}\otimes F_{3}\otimes L^{-1}M^{2} \cong 
(F_{2}\oplus 
F_{4})\otimes L^{-1}M^{2}$ (cf. \cite{At} p.437). It holds 
$h^0(\Check{F}\otimes S^2E_{2})\leq 2$, so one has only $\geq $ in 
\eqref{es2.17}. However here $\#\, \text{moduli}\, (E) \leq 2$, while 
$h^1(End\, E_{1})=1,\; h^1(End\, E_{2}) = h^0(End\, F_{2})=2$. Thus 
one has strict inequality in \eqref{es2.16} and \eqref{es2.12} holds.

{\bf Subcase 2C$'$(iii).}\quad $E_{2}\cong M_{1}\oplus M_{2}, \; M_{1}\ncong 
M_{2}$. Here $h^0(E_1\otimes E_2^*)\geq 2$ and $S^{2}E_{2}\cong 
M_{1}^{2}\oplus M_{1}M_{2}\oplus M_{2}^{2}$. Only two of these invertible 
sheaves might be isomorphic to each other, so $h^0(\Check{F}\otimes 
S^2E_{2})\leq 2$. Thus \eqref{es2.17} holds

{\bf Subcase 2C$'$(iv).}\quad $E_{2}\cong M\oplus M$. Here 
$h^0(\Check{F}\otimes S^2E_{2})\leq 3,\; h^1(End\, E_{1})=1, \; 
h^1(End\, E_{2})=4,\; h^0(E_1\otimes E_2^*)\geq 2$ and $\#\, 
\text{moduli}\, (E,F) = \#\, \text{moduli}\, (E) \leq 2$. Thus 
we have $\geq 0$ in \eqref{es2.17} and $< 0$ in \eqref{es2.16} so
the 
inequality \eqref{es2.12} holds.

\medskip
\noindent
Under the condition of Subcase~2C we now assume

{\bf Subcase 2C$''$.} \quad \emph{$F$ is decomposable, $F\cong L_{1}\oplus 
L_{2},\; \deg L_{i}=e/2$}. The inequality \eqref{es2.12} which we want 
to verify reads here as follows:
\begin{equation}\label{es2.18}
\begin{split}
\#\, \text{moduli}\, (E,F) 
&< h^1(End\, E_{1}) + h^1(End\, E_{2})\\
&\quad + h^1(End\, F)-1+h^0(E_1\otimes E_2^*)-h^0(\Check{F}\otimes 
S^2E_{2})
\end{split}
\end{equation}

{\bf Subcase 2C$''$(i).}\quad \emph{Let $r(E_{2})=1$}. Here $\det F\cong 
\det E$ yields $L_{1}\otimes L_{2}\cong \det E_{1}\otimes E_{2}$. One 
has $\Check{F}\otimes S^2E_{2}\cong (L_{1}^{-1}\otimes E_{2}^{2})\oplus 
(L_{2}^{-1}\otimes E_{2}^{2})$. So in order that $h^0(\Check{F}\otimes 
S^2E_{2})>0$ it should hold $L_{i}\cong E_{2}^{2}$ for $i=1$ or $i=2$. 
We conclude $F$ is determined uniquely by $E$, so $\#\, 
\text{moduli}\, (E,F) = \#\, \text{moduli}\, (E)$. We have 
\begin{equation}\label{es2.19}
\#\, \text{moduli}\, (E) \leq h^1(End\, E_{1}) + h^1(End\, E_{2})
\end{equation}
so proving that 
\begin{equation}\label{es2.19a}
h^0(End\, F) - 1 + h^0(E_1\otimes E_2^*) - h^0(\Check{F}\otimes 
S^2E_{2}) > 0
\end{equation}
would imply \eqref{es2.18}. We have 
\begin{equation}\label{es2.19b}
h^0(End\, F) - 1 = \begin{cases}
1& \text{if $F\cong L_{1}\oplus L_{2},\; L_{1}\ncong L_{2}$}\\
3& \text{if $F\cong L\oplus L$}
\end{cases}
\end{equation}
and $h^0(E_1\otimes E_2^*)=2(\mu_{1}-\mu_{2})\geq 1$. If 
$F\cong L_{1}\oplus L_{2},\; L_{1}\ncong L_{2}$, then 
$h^0(\Check{F}\otimes E_{2}^{2})\leq 1$ and if $F\cong L\oplus L$ then 
$h^0(\Check{F}\otimes E_{2}^{2})\leq 2$. In both cases \eqref{es2.19a} 
holds.

{\bf Subcase 2C$''$(ii).}\quad \emph{$r(E_{2})=2$}. Here $\Check{F}\otimes 
S^2E_{2}$ is a direct sum of indecomposable locally free sheaves of 
degree 0. An easy check of the various cases for $E_{2}$ shows that 
the conditions $h^0(\Check{F}\otimes S^2E_{2})>0$ and $\det F\cong 
\det E$ determine $F$ by $E$. We have again the inequality 
\eqref{es2.19} and it suffices to prove \eqref{es2.19a}. We proceed 
similarly to Subcases 2C$'$(i) -- 2C$'$(iv).

If $E_{2}$ is indecomposable of odd degree, then 
\(S^{2}E_{2}\cong \wedge^{2}E_{2}\otimes (\eta_{1}\oplus\eta_{2}\oplus 
\eta_{3})\) where $\eta_{i}$ are the three points of order 2 of 
$Pic^{0}Y$. If $F\cong L\oplus L$ then $h^0(\Check{F}\otimes 
S^2E_{2})\leq 2$. Since $h^0(End\, F)=4$ we obtain \eqref{es2.19a}. If 
$F\cong L_{1}\oplus L_{2}, \; L_{1}\ncong L_{2}$ then in case 
$h^0(\Check{F}\otimes S^2E_{2})\leq 1$ the inequality \eqref{es2.19a} 
holds. The only other possibility might be $h^0(\Check{F}\otimes 
S^2E_{2})=2$ when $L_{1}\cong \wedge^{2}E_{2}\otimes \eta_{i},\; 
L_{2}\cong \wedge^{2}E_{2}\otimes \eta_{j}$ for some pair $i,j$. Here 
we can only claim $\geq 0$ in \eqref{es2.19a}. Now the isomorphism 
$L_{1}L_{2}\cong \det F \cong \det E\cong E_{1}\otimes \det E_{2}$ 
yields $E_{1}\cong \wedge^2E_2\otimes \eta_i\eta_j$. Thus $\#\, \text{moduli}\, (E)\leq 1$. In this case 
the inequality \eqref{es2.19} is strict and again we see 
\eqref{es2.18} holds. 

If $E_{2}$ is indecomposable of even degree, $E_{2}\cong F_{2}\otimes 
M$, then $S^{2}E_{2}\cong F_{3}\otimes M^{2}$, so if $F\cong 
L_{1}\oplus L_{2}$ with $L_{1}\ncong L_{2}$ then $h^0(\Check{F}\otimes 
S^2E_{2})\leq 1$. If $F\cong L\oplus L$ then $h^0(\Check{F}\otimes 
S^2E_{2})\leq 2$. In both cases \eqref{es2.19a} holds.

If $E_{2}$ is decomposable, $E_{2}\cong M_{1}\oplus M_{2}$, then 
$h^0(E_1\otimes E_2^*)\geq 2,\; S^{2}E_{2}\cong M_1^2\oplus M_1M_2\oplus M_2^2$. If 
$h^0(\Check{F}\otimes S^2E_{2})\leq 2$ then \eqref{es2.19a} holds. One 
might have $h^0(\Check{F}\otimes S^2E_{2})\geq 3$ only if two of the 
summands of $S^{2}E_{2}$ are isomorphic to each other, i.e. either 
$M_{1}\ncong M_{2},\; M_{1}^{2}\cong M_{2}^{2}$ or $M_{1}\cong 
M_{2}\cong M$. In the former case one has $h^0(\Check{F}\otimes 
S^2E_{2})\geq 3$ either if $L_{1}\cong L_{2}\cong M_{1}^{2}\cong M_{2}^{2}$ 
and then $h^0(\Check{F}\otimes S^2E_{2})=4$ and \eqref{es2.19a} holds 
or if $L_{1}\cong M_{1}^{2}\cong M_{2}^{2},\; L_{2}\cong M_{1}M_{2}$ 
and then $h^0(\Check{F}\otimes S^2E_{2})=3$. In this case the 
left-hand side of \eqref{es2.19a} is $\geq 0$. Here one has 
$L_{1}L_{2}\cong \det E\cong E_{1}M_{1}M_{2}$.  Thus $E_{1}\cong M_{1}^{2}$ and $\#\, \text{moduli}\, (E)\leq 
1$ since $M_{1}^{2}\cong M_{2}^{2}$. We have that $h^1(End\, 
E_{1})+h^1(End\, E_{2})=3$, so \eqref{es2.18} holds. For $E_{2}$ it 
remains the case $E_{2}\cong M\oplus M$. One might have 
$h^0(\Check{F}\otimes S^2E_{2})\geq 3$ in two cases. Either $F\cong 
L_{1}\oplus L_{2},\; L_{1}\ncong L_{2},\; L_{1}\cong M^2$ and then 
$h^0(\Check{F}\otimes S^2E_{2})=3$, or $F\cong L\oplus L,\; L\cong M^2$ 
and then $h^0(\Check{F}\otimes S^2E_{2})=6$. In both cases the 
left-hand side of \eqref{es2.19a} is $\geq -1$ while $\#\, 
\text{moduli}\, (E,F)=\#\, \text{moduli}\, (E)\leq 2$ and $h^1(End\, 
E_{1})=1,\; h^1(End\, E_{2})=4$. Thus \eqref{es2.18} holds.
Case 2 is completed.
\end{block}

\begin{block}\label{s2.21}
{\bf Case 3.} \emph{$F$ is semistable, $E=E_{1}\oplus E_{2}\oplus E_{3}$ 
with $d(E_{1})>d(E_{2})>d(E_{3})$}. Let $d_{i}=d(E_{i})$. We have 
$h^1(End\, E) = \sum h^1(End\, E_{i}) + \sum_{i<j}h^0(E_i\otimes 
E_j^*) = 3 +2(d_{1}-d_{3})$, so the inequality \eqref{es2.12} that we 
want to verify for different types of $(E,F)$ within this case becomes 
\begin{equation}\label{es2.21}
\#\, \text{moduli}\, (E,F)\ <\ 2+h^1(End\, 
F)+2(d_{1}-d_{3})-h^1(\Check{F}\otimes S^2E)
\end{equation}

{\bf Subcase 3A.}\quad \emph{$\mu(\Check{F}\otimes E_{3}^{2})= 
-\frac{e}{2}+2d_{3}>0$ or $\mu(\Check{F}\otimes E_{3}^{2})=0$ and 
$h^0(\Check{F}\otimes E_{3}^{2})=0$}. Here we have 
$h^1(\Check{F}\otimes S^2E)=0,\; \#\, \text{moduli}\, (E,F)\leq 
3+h^1(End\, F)-1$  and $d_{1}-d_{3}\geq 2$, thus \eqref{es2.21} holds.

{\bf Subcase 3B.}\quad $\mu(\Check{F}\otimes E_{3}^{2})<0$. The same 
argument as in Subcase 2B of (\ref{s2.15}) shows that such types of 
$(E,F)$ cannot occur in the case of pair associated with an 
irreducible, Gorenstein quadruple cover $X$.

{\bf Subcase 3C.}\quad $\mu(\Check{F}\otimes E_{3}^{2})=0,\; 
h^0(\Check{F}\otimes E_{3}^{2})>0$. Here one has \linebreak
$h^1(\Check{F}\otimes 
S^2E) = 
h^1(\Check{F}\otimes E_{3}^{2}) = h^0(\Check{F}\otimes 
E_{3}^{2})\leq 2$. Since $d_{1}-d_{3}\geq 2$ we conclude 
\eqref{es2.21} holds.
\end{block}

\bigskip
\noindent
We thus examined in (\ref{s2.12}) -- (\ref{s2.21}) al cases for 
$(E,F)$ with semistable $F$. The remaining types to be considered are 
with $F\cong L_{1}\oplus L_{2}$ where $d(L_{1})>d(L_{2})$. Let 
$d(L_{i})=\lambda_{i}$. Here $h^1(End\, F)=2+h^0(L_{1}L_{2}^{-1}) = 2 
+ \lambda_{1}-\lambda_{2}$. A possible type for $(E,F)$ would yield 
insufficient number of moduli if it holds the inequality (cf. 
\eqref{es2.12})
\begin{equation}\label{es2.22}
\#\, \text{moduli}\, (E,F)\ <\ h^1(End\, 
E)+1+(\lambda_{1}-\lambda_{2})- h^1(\Check{F}\otimes S^2E)
\end{equation}
It obviously holds $\#\, \text{moduli}\, (E,F) \leq h^1(End\, E) + 1$

\begin{block}\label{s2.22}
{\bf Case 4.}\quad \emph{$F\cong L_{1}\oplus L_{2}$, $E$ is semistable}. Since 
$\det E\cong \det F$ we have $e=\lambda_{1}+\lambda_{2}=3\mu(E)$ and 
$\lambda_{1}>\frac{e}{2}>\lambda_{2}$. One has $h^1(\Check{F}\otimes 
S^2E) = h^1(L_{1}^{-1}\otimes S^{2}E)+h^1(L_{2}^{-1}\otimes S^2E)$ and 
$\mu(L_{2}^{-1}\otimes S^2E)=-\lambda_{2}+\frac{2e}{3}>\frac{e}{6}>0$. 
Thus $h^1(L_{2}^{-1}\otimes S^2E)=0$ and the inequality \eqref{es2.22} 
we aim to prove becomes
\begin{equation}\label{es2.22a}
\#\, \text{moduli}\, (E,F) \ <\ h^1(End\, E) + 1 + 
(\lambda_{1}-\lambda_{2}) - h^1(L_{1}^{-1}\otimes S^2E)
\end{equation}

{\bf Subcase 4A.} \quad $\mu(L_{1}^{-1}\otimes S^2E)>0$. Here $h^1(
L_{1}^{-1}\otimes S^2E
)=0$, therefore $h^1(\Check{F}\otimes S^2E)=0$ and \eqref{es2.22a} 
holds.

{\bf Subcase 4B.}\quad \emph{$\mu(L_{1}^{-1}\otimes S^2E)<0$ or 
$\mu(L_{1}^{-1}\otimes S^2E)=0$ and $h^0(L_{1}^{-1}\otimes S^2E)=0$.}
 Here we have $Hom_{Y}(L_{1},S^{2}E)=0$. This is impossible for a 
pair associated to a Gorenstein quadruple covering of $Y$. Indeed, the 
associated relative pencil of conics is given by a monomorphism $\eta: 
F\to S^{2}E$ (cf. \cite{CE} pp.450,451), so $\eta|_{L_{1}}$ must be a nonzero 
element for such a pair.

{\bf Subcase 4C.}\quad \emph{$\mu(L_{1}^{-1}\otimes S^2E)=0$ and 
$h^0(L_{1}^{-1}\otimes S^2E)>0$}. Here $d(L_{1}^{-1}\otimes 
S^{2}E)=0$, so $h^1(L_{1}^{-1}\otimes S^2E)=h^0(L_{1}^{-1}\otimes 
S^2E)$. Furthermore $\lambda_{1}=\mu(L_{1})=\mu(S^{2}E)=\frac{2e}{3}$, 
so in this subcase $e\equiv 0(mod\, 3)$. 

{\bf Subcase 4C(i).}\quad \emph{$E$ is indecomposable}. Here $E\cong 
F_{3}\otimes M$ since $3|e$. We have $S^{2}E=S^{2}F_{3}\otimes M^{2}$ 
and we claim $S^{2}F_{3}\cong F_{1}\oplus F_{5}$. Indeed, 
$F_{3}\otimes F_{3}\cong F_{1}\oplus F_{3}\oplus F_{5}$ by 
\cite{At} p.438. One has $F_{3}\otimes F_{3}\cong S^{2}F_{3}\oplus 
\wedge^{2}F_{3}$. The pairing $\wedge^{2}F_{3}\times F_{3}\to 
\wedge^{3}F_{3}\cong \mathcal{O}_{Y}$ is nondegenerate, so 
$\wedge^{2}F_{3}\cong F_{3}^{*}\cong F_{3}$. We conclude by  \cite{At1} that $S^{2}F_{3}\cong F_{1}\oplus 
F_{5}$. Calculating the various entries in \eqref{es2.22a} we obtain: 
$h^0(L_{1}^{-1}\otimes S^2E)\leq 2, \; h^1(End\, E) = 3$, so the 
right-hand side is $\geq 3$. The condition $\det E\cong \det F$ reads 
as $M^{3}\cong L_{1}L_{2}$, so $E$ is determined by $F$ up to 
tensoring by a ponint of order 3 in $Pic^{0}Y$. Hence $\#\, 
\text{moduli}\, (E,F) = \#\, \text{moduli}\, (F) \leq 2$ and the 
inequality \eqref{es2.22} holds.

{\bf Subcase 4C(ii).}\quad \emph{$E\cong E_{1}\oplus E_{2}$ where $E_{1}$ is 
indecomposable of rank 2}. We have from semistability 
$\mu(E)=\mu(E_{1})=\mu(E_{2})$, thus $\mu(E_{1})=\frac{e}{3}$. Since 
$3|e$ we conclude $d(E_1)$ is even, so $E_{1}\cong F_{2}\otimes M_{1}$. 
We have $h^1(End\, E)=h^0(End\, F_{2})+h^0(End\, 
E_{2})+h^0(F_{2}\otimes M_{1}E_{2}^{-1})$, so $h^1(End\, E) = 3$ if 
$M_{1}\ncong E_{2}$ and $h^1(End\, E) = 4$ if $M_{1}\cong E_{2}$. We 
have $S^{2}E\cong \linebreak
S^{2}E_{1}\oplus (E_{1}\otimes E_{2})\oplus 
E_{2}^{2}$, so 
\[
L_{1}^{-1}\otimes S^{2}E \cong (F_{3}\otimes L_{1}^{-1}M_{1}^{2})\oplus 
(F_{2}\otimes L_{1}^{-1}M_{1}E_{2})\oplus L_{1}^{-1}E_{2}^{2}
\]
We have now various cases:

If $M_{1}^{2}, M_{1}E_{2}$ and $E_{2}^{2}$ are not isomorphic to each 
other then $h^0(L_{1}^{-1}\otimes S^{2}E )\leq 1$. So the right-hand 
side of \eqref{es2.22a} is $\geq 4$ while $\#\, \text{moduli}\, 
(E,F)\leq 2 + 2 -1 =3$, thus \eqref{es2.22a} holds.

If $M_{1}^{2}\cong E_{2}^{2}$ but $M_{1}\ncong E_{2}$ then 
$h^0(L_{1}^{-1}\otimes S^{2}E )\leq 2$, while $\#\, \text{moduli}\, (E,F) \linebreak
\leq 
2$, so again \eqref{es2.22a} holds. 

If $M_{1}\cong E_{2}$, then $h^0(L_{1}^{-1}\otimes S^{2}E )\leq 3,\; h^1(End\, 
E) = 4$ and $\#\, \text{moduli}\, (E,F) \linebreak
\leq 2$, so \eqref{es2.22a} 
holds.

{\bf Subcase 4C(iii).}\quad \emph{$E\cong M_{1}\oplus M_{2}\oplus M_{3}$ 
where $\deg M_{i} = \frac{e}{3}$}. Here one has $h^1(End\, E) = 3 + 2 
\sum_{i<j}h^0(M_{i}M_{j}^{-1})$ and $S^{2}E\cong \oplus _{i\leq 
j}M_{i}M_{j}$. The conditions $h^0(L_{1}^{-1}\otimes S^{2}E)\geq 1$ and $\det 
F\cong \det E$ determine $F$ from $E$. Hence $\#\, \text{moduli}\, 
(E,F) = \#\, \text{moduli}\, (E)\leq 3\leq h^1(End\, E)$. If $h^0(
L_{1}^{-1}\otimes S^{2}E) = 1$ we see \eqref{es2.22a} holds. One has the 
inequality $h^0(L_{1}^{-1}\otimes S^{2}E)\geq 2$ in the following cases.

a) Up to reordering $M_{1}\cong M_{2}, M_{3}\ncong M_{1}$. Then 
$S^2E\cong (M_{1}^2)^{\oplus 3}\oplus 
(M_{1}M_{3})^{\oplus 2}\oplus M_{3}^{2}$. Here 
$h^0(L_{1}^{-1}\otimes S^{2}E)\leq 3, h^1(End\, E) = 5$, so \eqref{es2.22a} holds.

b) The sheaves $M_{1},M_{2},M_{3}$ are pairwise non-isomorphic, two of 
the summands of $S^2E$ are isomorphic to each other, $L_{1}$ is 
isomorphic to this summand and no three summands of $S^2E$ are 
isomorphic to each other. Up to reordering this might happen if 
$M_{1}^{2}\cong M_{2}^{2}\cong L_{1}$ or $M_{1}^2\cong M_{2}M_{3}\cong 
L_{1}$. Then $h^0(L_{1}^{-1}\otimes S^{2}E)=2,\; \#\, \text{moduli}\, 
(E)\leq 2$ thus \eqref{es2.22a} holds as well.

c) The sheaves $M_{1},M_{2},M_{3}$ are pairwise non-isomorphic, 
$M_{1}^{2}\cong M_{2}^{2}\cong M_{3}^{2}\cong L_{1}$. 
Here $h^0(L_{1}^{-1}\otimes S^{2}E)=3,\; \#\, \text{moduli}\, (E)\leq 
1$, thus \eqref{es2.22a} holds.

d) $M_{1}\cong M_{2}\cong M_{3},\; L_{1}\cong M_{1}^{2}$. Here 
$h^0(L_{1}^{-1}\otimes S^2E)=6,\; h^1(End\, E)=9$ and $\#\, 
\text{moduli}\, (E)\leq 1$ thus \eqref{es2.22a} holds as well. Case 4 
is completed.
\end{block}

\begin{block}\label{s2.26}
{\bf Case 5.}\quad \emph{$F\cong L_{1}\oplus L_{2},\; E\cong E_{1}\oplus 
E_{2}$ where $E_{1}$ is semistable of rank 2 and 
$\mu(E_{1})=\mu_{1}>\mu_{2}=d(E_{2})$}. Here $h^1(End\, E)=h^1(End\, 
E_{1})+1+2(\mu_{1}-\mu_{2})$. So, a possible type of $(E,F)$ of this 
class  would be excluded if one shows that 
\begin{equation}\label{es2.26}
\begin{split}
\#\, \text{moduli}\, (E,F) < & h^1(End\, E_{1})+2\\
&+(\lambda_{1}-\lambda_{2}) + 2(\mu_{1}-\mu_{2}) - 
h^1(\Check{F}\otimes S^2E).
\end{split}
\end{equation}
The condition $\det E\cong \det E_{1}\otimes E_{2}\cong \det F\cong 
L_{1}L_{2}$ determines $E_{2}$ from $E_{1}$ and $F$, so 
\begin{equation}\label{es2.26a}
\#\, \text{moduli}\, (E,F) \ \leq \ h^1(End\, E_{1}) + 2
\end{equation}
We see \eqref{es2.26} would hold in cases when
\begin{equation}\label{es2.26b}
\lambda_{1}-\lambda_{2}+2(\mu_{1}-\mu_{2}) - h^1(\Check{F}\otimes 
S^2E) > 0.
\end{equation}
One has 

\begin{equation*}
\begin{split}
h^1(\Check{F}\otimes S^2E)
&=h^1(L_{1}^{-1}\otimes S^2E_{1})+
h^1(L_{1}^{-1}\otimes E_{1}\otimes E_{2})+
h^1(L_{1}^{-1}\otimes E_{2}^{2})\\
&+h^1(L_{2}^{-1}\otimes S^2E_{1})+ h^1(L_{2}^{-1}\otimes 
E_{1}\otimes E_{2})+h^1(L_{2}^{-1}\otimes E_{2}^{2}).
\end{split}
\end{equation*}
The direct summands $L_{1}^{-1}\otimes S^2E_{1},\ldots,
L_{2}^{-1}\otimes E_{2}^{2}$ are semistable with the following slopes
\begin{equation}\label{es2.26c}
\begin{split}
&-\lambda_{1}+2\mu_{1} > -\lambda_{1}+\mu_{1}+\mu_{2} > 
-\lambda_{1}+2\mu_{2}\\
&-\lambda_{2}+2\mu_{1} > -\lambda_{2}+\mu_{1}+\mu_{2} > 
-\lambda_{2}+2\mu_{2}
\end{split}
\end{equation}
where in each column the upper number is smaller then the lower one. 
If $-\lambda_{2}+2\mu_{2}<0$ or if $-\lambda_{2}+2\mu_{2}=0$ and 
$h^0(L_{2}^{-1}\otimes E_{2}^{2})=0$ then $h^0(\Check{F}\otimes 
E^{2}_{2})$ = 0 and the same argument as in Subcase 2B of 
(\ref{s2.15}) shows that such case is impossible for a pair obtained 
from an irreducible cover, so it is excluded. 
If $-\lambda_{1}+\mu_{1}+\mu_{2}<0$ or if 
$-\lambda_{1}+\mu_{1}+\mu_{2}=0$ and 
$h^1(L_{1}^{-1}\otimes E_{1}\otimes E_{2})=0$ then 
$H^0(L_{1}^{-1}\otimes E_{1}\otimes E_{2})=0=
H^0(L_{1}^{-1}\otimes E_{2}^{2})$ or equivalently $H^0(L_{1}^{-1}E_{2}\otimes E)=0$. 
Suppose that in the 
explicit representation of  (\ref{s2.14}) the sections 
$e_{1},e_{2}$ form a frame of $E_{1}|_{U}$ and $e_{3}$ generates 
$E_{2}|_{U}$. Then representing $\eta$ by \eqref{es2.14} we see that for 
each $y\in U$ the quadratic polinomial 
$\sum_{i\leq j}a_{ij,1}e_{i}(y)e_{j}(y)$ has no monomials which 
contain $e_{3}(y)$. Hence the corresonding conic is degenerate. If 
$\eta \in H^{0}(\Check{F}\otimes S^2E) = Hom_{Y}(L_{1}\oplus L_{2}, 
S^{2}E)$ is obtained from a reduced cover $X$, then $\eta|_{L_{1}}$ 
yields a conic bundle over $Y$ whose general fibre is a union of two 
distinct lines. Let $\pi_{*}:\tilde{Y}\to Y$ be the associated double 
covering. We obtain that $\pi :X\to Y$ may be decomposed as $X\to 
\tilde{Y}\to Y$. We consider only simple coverings, so only the case 
of \'{e}tale $\pi_{2}:\tilde{Y}\to Y$ is of interest to us. Now, such 
coverings are excluded from the hypothesis of Theorem~\ref{s2.9}, 
namely from the condition that $\pi_*:H_1(X,\mathbb{Z})\to 
H_1(Y,\mathbb{Z})$ is surjective. We conclude that pairs $(E,F)$ of 
Case 5 may yield quadruple coverings $\pi :X\to Y$ from 
$\mathcal{H}^{0}_{4,n}(Y)$ only if 
\begin{itemize}
\item $-\lambda_{2}+2\mu_{2}\geq 0$ and if it is 0 
then $h^0(L_{2}^{-1}\otimes E_{2}^{2})\geq 1$.
\item $-\lambda_{1}+\mu_{1}+\mu_{2}\geq 0$ and if it is 0 then
$h^0(L_{1}^{-1}\otimes E_{1}\otimes E_{2})\geq 1$.
\end{itemize}
Looking at \eqref{es2.26c} we obtain
\begin{equation}\label{es2.28}
h^1(\Check{F}\otimes S^2E)
=
h^1(L_{1}^{-1}\otimes E_{1}\otimes E_{2})+
h^1(L_{1}^{-1}\otimes E_{2}^{2})
+h^1(L_{2}^{-1}\otimes E_{2}^{2}).
\end{equation}
Notice that 
\(-\lambda_{2}+2\mu_{1}>-\lambda_{1}+2\mu_{1}>
 -\lambda_{1}+\mu_{1}+\mu_{2}\geq 0\). Since $2\mu_{1}\in \mathbb{Z}$ 
we conclude $-\lambda_{2}+2\mu_{1}\geq 2$.

{\bf Subcase 5A.}\quad \emph{Assume $-\lambda_{2}+2\mu_{2}>0$ and 
$-\lambda_{1}+\mu_{1} + \mu_{2}>0$}. Then 
$h^1(\Check{F}\otimes S^2E)=h^1(L_{1}^{-1}\otimes E_{2}^{2})$. The 
inequality \eqref{es2.26b} holds obviously if 
$-\lambda_{1}+2\mu_{2}\geq 0$. If $-\lambda_{1}+2\mu_{2}< 0$ then 
$h^1(L_{1}^{-1}\otimes E_{2}^{2})=\lambda_{1}-2\mu_{2}$ and we have
\begin{equation*}
\lambda_{1}-\lambda_{2} + 2(\mu_{1}-\mu_{2}) - 
(\lambda_{1}-2\mu_{2})=-\lambda_{2}+2\mu_{1}\geq 2
\end{equation*}
Therefore \eqref{es2.26b} holds.

{\bf Subcase 5B.}\quad \emph{Assume $-\lambda_{1}+\mu_{1}+\mu_{2}>0,\; 
-\lambda_{2}+2\mu_{2}=0$ and $h^0(L_{2}^{-1}\otimes E_{2}^{2})\geq 
1$}. Then $h^1(\Check{F}\otimes S^2E)=\lambda_{1}-2\mu_{2}+1$. As in 
the preceeding case we obtain \eqref{es2.26b} holds.

{\bf Subcase 5C.}\quad \emph{Assume $-\lambda_{2}+2\mu_{2}>0, \; 
-\lambda_{1}+\mu_{1}+\mu_{2}=0$ and 
$h^0(L_{1}^{-1}\otimes E_{1}\otimes E_{2})\geq 1$}. Then 
$h^1(\Check{F}\otimes S^2E)=
h^0(L_{1}^{-1}\otimes E_{1}\otimes E_{2})+\lambda_{1}-2\mu_{2}$. If 
$h^0(L_{1}^{-1}\otimes E_{1}\otimes E_{2})\leq 1$ we conclude as in 
the preceeding case. The only other possibility might be $E_{1}\cong 
M\oplus M$ and $h^0(L_{1}^{-1}\otimes E_{1}\otimes E_{2})=2$. However 
then the inequality \eqref{es2.26a} is strict since $h^1(End\, 
E_{1})=4$. Thus \eqref{es2.26} holds.

{\bf Subcase 5D.}\quad \emph{Assume $-\lambda_{2}+2\mu_{2}=0,\; 
h^0(L_{2}^{-1}\otimes E_{2}^{2})\geq 1$ and 
$-\lambda_{1}+\mu_{1}+\mu_{2}=0, 
h^0(L_{1}^{-1}\otimes E_{1}\otimes E_{2})\geq 1$}. Then 
\[
h^1(\Check{F}\otimes S^2E)
=
h^0(L_{1}^{-1}\otimes E_{1}\otimes E_{2})+
h^0(L_{1}^{-1}\otimes E_{2}^{2}) + \lambda_{1}-2\mu_{2}
\]
The right-hand side of \eqref{es2.26} becomes
\[
h^1(End\, E_{1}) + 2 + (-\lambda_{2}+2\mu_{1}) - 
(h^0(L_{1}^{-1}\otimes E_{1}\otimes E_{2})+
h^0(L_{1}^{-1}\otimes E_{2}^{2}))
\]
As we saw $-\lambda_{2}+2\mu_{1}\geq 2$, while 
$2\leq 
h^0(L_{1}^{-1}\otimes E_{1}\otimes E_{2})+
h^0(L_{1}^{-1}\otimes E_{2}^{2})\leq 3
$.

{\bf Subcase 5D(i)}\quad $h^0(L_{1}^{-1}\otimes E_{2}^{2})=1,\;
h^0(L_{1}^{-1}\otimes E_{1}\otimes E_{2})=1$. Here from the first 
equality $L_{2}\cong E_{2}^{2}$ and from $\det F\cong \det E$ it 
follows $L_{2}\cong \det E_{1}\otimes E_{2}^{-1}$. Thus $\#\, 
\text{moduli}\, (E,F)\leq h^1(End\, E_{1})+1$ and therefore 
\eqref{es2.26} holds.

{\bf Subcase 5D(ii)}\quad 
$h^0(L_{1}^{-1}\otimes E_{2}^{2})=1,\;
h^0(L_{1}^{-1}\otimes E_{1}\otimes E_{2})=2$. With respect to the 
previuos case we have in addition that $E_{1}\cong M\oplus M$, so 
$h^1(End\, E_{1})=4$ and \eqref{es2.26} is fulfilled. Case 5 is 
completed.
\end{block}

\begin{block}\label{s2.30}
{\bf Case 6.}\quad \emph{$F\cong L_{1}\oplus L_{2},\; E\cong E_{1}\oplus 
E_{2}$ where $E_{2}$ is semistable of rank 2 and 
$\mu_{1}=d(E_{1})>\mu_{2}=\mu(E_{2})$}. A calculation similar to the 
one in (\ref{s2.26}) shows that a type of $(E,F)$ in this class would 
be excluded if one shows that 
\begin{equation}\label{es2.30}
\begin{split}
\#\, \text{moduli}\, (E,F) < &h^1(End\, E_{2}) + 2\\ 
&+(\lambda_{1}-\lambda_{2})+2(\mu_{1}-\mu_{2})-h^1(\Check{F}\otimes 
S^2E)
\end{split}
\end{equation}
Furthermore it holds $\#\, \text{moduli}\, (E,F)\leq h^1(End\, 
E_{2})+2$, so \eqref{es2.30} would hold if one shows that
\begin{equation}\label{es2.30a}
(\lambda_{1}-\lambda_{2})+2(\mu_{1}-\mu_{2})-h^1(\Check{F}\otimes 
S^2E) > 0
\end{equation}
One has
\begin{equation*}
\begin{split}
h^1(\Check{F}\otimes S^2E)
&=h^1(L_{1}^{-1}\otimes E_{1}^{2})+
h^1(L_{1}^{-1}\otimes E_{1}\otimes E_{2})+
h^1(L_{1}^{-1}\otimes S^2E_{2})\\
&+h^1(L_{2}^{-1}\otimes E_{1}^{2})+ 
h^1(L_{2}^{-1}\otimes E_{1}\otimes E_{2})
+h^1(L_{2}^{-1}\otimes S^2E_{2})
\end{split}
\end{equation*}
where for the slopes of the semistable direct summands of 
$\Check{F}\otimes S^2E$ one has the same table as that of 
\eqref{es2.26c}. If $-\lambda_{1}+2\mu_{2}<0$ or if 
$-\lambda_{1}+2\mu_{2}=0$ and $h^0(L_{1}^{-1}\otimes S^2E_{2})=0$ then 
one would have $H^0(L_{1}^{-1}\otimes S^2E_{2})=0$. 
We may 
choose 
in (\ref{s2.14})
 $e_{1}$ to generate $E_{1}|_{U}$ and $e_{1},e_{2}$ to form a 
frame of $E_{2}|_{U}$. Then every 
 $\eta \in 
H^0(\Check{F}\otimes S^2E)=Hom_{Y}(F,S^2E)$ has the property 
that in the explicit representation of \eqref{es2.14} all 
monomials of $\eta|_{L_{1}}$ contain $e_{1}(y)$ as a factor. This 
means that the conic bundle corresponding to $\eta|_{L_{1}}$ is 
reducible and one of its components is the ruled surface corresponding 
to the epimorphism $E\to E_{1}\to 0$. This implies that the quadruple 
cover $X$ is reducible, which is excluded from the hypothesis of 
Theorem~\ref{s2.9}. 

We may thus assume that $-\lambda_{1}+2\mu_{2}\geq 
0$ and if $-\lambda_{1}+2\mu_{2}=0$ then 
$h^0(L_{1}^{-1}\otimes S^2E_{2})\geq 1$. Looking at the table 
\eqref{es2.26c}
 we see that $h^1(\Check{F}\otimes S^2E)=
h^1(L_{1}^{-1}\otimes S^2E_{2})$. If $-\lambda_{1}+2\mu_{2}>0$ then $
h^1(L_{1}^{-1}\otimes S^2E_{2})=0$,
 so \eqref{es2.30a} holds. Assume $-\lambda_{1}+2\mu_{2}=0$. Then 
$h^1(L_{1}^{-1}\otimes S^2E_{2})=h^0(L_{1}^{-1}\otimes S^2E_{2})$. Since \linebreak
$L_{1}^{-1}\otimes S^2E_{2}$ is semistable of rank 3 and slope 0 one 
has $1\leq \linebreak
h^0(L_{1}^{-1}\otimes S^2E_{2})\leq 3$. If 
$h^0(L_{1}^{-1}\otimes S^2E_{2})=1$ then \eqref{es2.30a} holds since 
$(\lambda_{1}-\lambda_{2})+2(\mu_{1}-\mu_{2})\geq 2$. This case 
happens if $E_{2}$ is indecomposable. Indeed, in this case if 
$d(E_{2})$ is even, then $S^{2}E_{2}\cong F_{3}\otimes M^{2}$, so by 
\cite{At}\quad $h^0(F_{3}\otimes M^{2}L_{1}^{-1})\leq 1$. If $d(E_{2})$ 
is odd then by Proposition~\ref{s1.5a} one has $S^2E_{2}\cong 
\wedge^{2}E_{2}\otimes (\eta_{1}\oplus \eta_{2}\oplus \eta_{3})$ where 
$\eta_{i}$ are the three points of order 2 in $Pic^{0}Y$, so again 
$h^0(L_{1}^{-1}\otimes S^2E_{2})\leq 1$. So, cases with $
h^0(L_{1}^{-1}\otimes S^2E_{2})\geq 2
$ could occur only if $E_{2}$ is decomposable. Then 
$\mu_{2}=\mu(E_{2})$ is an integer, $\mu_{1}>\mu_{2}$, so 
$(\lambda_{1}-\lambda_{2})+2(\mu_{1}-\mu_{2})\geq 3$. If $E_{2}$ is 
decomposable and $h^0(L_{1}^{-1}\otimes S^2E_{2})\leq 2$ the 
inequality \eqref{es2.30a} holds. The only remaining subcase is 
$h^0(L_{1}^{-1}\otimes S^2E_{2})=3$ which is possible if $E_{2}\cong 
M\oplus M,\; L_{1}\cong M^{2}$. Then $\#\, \text{moduli}\, (E,F)\leq 
2$ while $h^1(End\, E_{2})=4$. Thus \eqref{es2.30} holds. All possible 
types of Case 6 are thus excluded.
\end{block}

\begin{block}\label{s2.32}
{\bf Case 7}.\quad \emph{$F\cong L_{1}\oplus L_{2},\; E\cong 
E_{1}\oplus E_{2}\oplus E_{3}$ where $d(E_{i})=\mu_{i},\; 
\mu_{1}>\mu_{2}>\mu_{3}$}. Here $h^1(End\, 
F)=2+\lambda_{1}-\lambda_{2},\; h^1(End\, E)=3+2(\mu_{1}-\mu_{3})$. 
So, a possible type for $(E,F)$ would be excluded if one shows that 
(cf. \eqref{es2.12})
\begin{equation}\label{es2.32}
\#\, \text{moduli}\, (E,F) < 4 + (\lambda_{1}-\lambda_{2})+
2(\mu_{1}-\mu_{3}) - h^1(\Check{F}\otimes S^2E)
\end{equation}
Since $\#\, \text{moduli}\, (E,F)\leq 4$ this would be the case if the following inequality 
holds. 
\begin{equation}\label{es2.32a}
(\lambda_{1}-\lambda_{2})+2(\mu_{1}-\mu_{3})-
h^1(\Check{F}\otimes S^2E) > 0
\end{equation}
We have $h^1(\Check{F}\otimes S^2E)=
\sum_{i=1}^{2}\sum_{j\leq k}h^1(L_{i}^{-1}E_{j}E_{k})$ and the direct 
summands of $\Check{F}\otimes S^2E$ have degrees
\begin{equation*}
\begin{split}
&-\lambda_{1}+2\mu_{1}>-\lambda_{1}+\mu_{1}+\mu_{2} >d_1>
-\lambda_{1}+\mu_{2}+\mu_{3}>-\lambda_{1}+2\mu_{3}\\
&-\lambda_{2}+2\mu_{1}>-\lambda_{2}+\mu_{1}+\mu_{2} >d_2>
-\lambda_{2}+\mu_{2}+\mu_{3}>-\lambda_{2}+2\mu_{3}.
\end{split}
\end{equation*}
Here $d_i$ stays for either $-\lambda_i + 2\mu_{2}$ or $-\lambda_{i}+\mu_{1}+\mu_{3}$ 
and in each column the upper number is smaller then the respective lower one.
If $-\lambda_{2}+2\mu_{3}<0$  then $H^0(\Check{F}\otimes E_{3}^2)=0$. If
$-\lambda_{1}+2\mu_{2}<0$  then $H^{0}(L_{1}^{-1}\otimes 
S^{2}(E_{2}\oplus E_{3}))=0$. 
If $-\lambda_{1}+\mu_{1}+\mu_{3}<0$ then $H^{0}(L_{1}^{-1}E_3\otimes E)=0$.
In each case the same arguments as 
those in Subcase 2B, Case 6 or Case 5 respectively show that such possibilities are 
excluded from the assumption that the quadruple cover is irreducible and cannot be decomposed through a double covering.
So we may assume $-\lambda_{2}+2\mu_{3}\geq 0$, 
$-\lambda_{1}+2\mu_{2}\geq 0$  and $-\lambda_{1}+\mu_{1}+\mu_{3}\geq 0$. We thus 
have
\begin{multline}\label{es2.33}
h^1(\Check{F}\otimes S^2E) = h^1(L_{1}^{-1}E_{1}E_{3})+
h^1(L_{1}^{-1}E_{2}^{2})\\ 
+h^1(L_{1}^{-1}E_{2}E_{3})+
h^1(L_{1}^{-1}E_{3}^{2})+ h^1(L_{2}^{-1}E_{3}^{2})
\end{multline}
We now consider various cases. If $-\lambda_{1}+2\mu_{3}\geq 0$ then 
all summands in \eqref{es2.33} are zero except possibly $
h^1(L_{1}^{-1}E_{3}^2)
$ which is $\leq 1$. So \eqref{es2.32a} holds. If 
$-\lambda_{1}+2\mu_{3}<0$ and $-\lambda_{1}+\mu_{2}+\mu_{3}\geq 0$ 
then $h^1(\Check{F}\otimes S^2E) = \epsilon_{1}+(\lambda_{1}-2\mu_{3})
+\epsilon_{2}$ where $0\leq \epsilon_{i}\leq 1$. The left-hand side of 
\eqref{es2.32a} becomes 
$-\lambda_{2}+2\mu_{3}+2(\mu_{1}-\mu_{3})-\epsilon_{1}-\epsilon_{2}
\geq 2$, so \eqref{es2.32a} holds. If 
$-\lambda_{1}+\mu_{2}+\mu_{3}<0$ then $h^1(\Check{F}\otimes S^2E) = 
\epsilon_{1}+\epsilon_{2}+(\lambda_{1}-\mu_{2}-\mu_{3})+(\lambda_{1}-2\mu_{3})+
\epsilon_{3}$ where $0\leq \epsilon_{i}\leq 1$. The left-hand side of 
\eqref{es2.32a} becomes 
\[
(-\lambda_{2}+2\mu_{3})+(-\lambda_{1}+2\mu_{1})+(\mu_{2}-\mu_{3})
-\epsilon_{1}-\epsilon_{2}-\epsilon_{3}
\]
We have $-\lambda_{2}+2\mu_{3}\geq 0,\; -\lambda_{1}+2\mu_{1}>
-\lambda_{1}+\mu_{1}+\mu_{2}>-\lambda_{1}+2\mu_{2}\geq 0$, so 
$-\lambda_{1}+2\mu_{1}\geq 2$. Moreover $\mu_{2}-\mu_{3}\geq 1$. Hence 
inequality \eqref{es2.32a} holds when at least one of $\epsilon_{i}$ is zero. If $\epsilon_{i}=1$ for $\forall i$ then we would have $L_1\cong E_1^2, L_1\cong E_1E_3, L_2\cong E_3^2$. In this case $\#\, \text{moduli}\, (E,F) \leq 2$ therefore \eqref{es2.32} holds. All possible types of Case 7 are 
thus excluded. 
\end{block}

\bigskip
\noindent
We thus verified that except for the types of $(E,F)$ specified in the the theorem and considered in Case 1 all other types yield number of parameters for the equivalence classes of coverings $X\to Y$ less then $n$. In order to complete the proof of the theorem
we need a lemma which is related to Theorem 4.5 of \cite{CE} and is an 
analog of Lemma 2.8 of \cite{K1}. We state and prove a more general 
result than we actually need for the proof of Theorem~\ref{s2.9}.

\begin{lem}\label{s2.34a}
Assume the base field $k$ is algebraically closed and $char(k)=0$. Let 
$q:\mathcal{Y}\to Z$ be a smooth, proper morphism with connected 
fibers, where $Z$ is smooth. Let $\mathcal{E}$ and 
$\mathcal{F}$ be locally free sheaves on $\mathcal{Y}$ of ranks 3 and 2 
respectively, such that $\det \mathcal{E}_{z}\cong \det \mathcal{F}_{z}$ 
for every $z\in Z$. Suppose 
$h^0(\mathcal{Y}_{z},\Check{\mathcal{F}}_{z}\otimes S^{2}\mathcal{E}_{z})$ 
is independent of $z\in Z$ and is $\ne 0$. 
Consider the locally free sheaf $\mathcal{H}=q_{*}
(\Check{\mathcal{F}}\otimes S^{2}\mathcal{E})$
 on $Z$. Let $f:\mathbb{H}\to Z$ be the associated 
vector bundle with fibers 
$\mathbb{H}_z=H^0(\mathcal{Y}_{z},
\Check{\mathcal{F}}_{z}\otimes S^{2}\mathcal{E}_{z}
)$. Then the subset $\mathbb{H}_{0}\subset \mathbb{H}$ 
consisting of $\eta$ which satisfy the following three conditions is 
Zariski open   in $\mathbb{H}$.
\renewcommand{\theenumi}{\alph{enumi}}
\begin{enumerate}
\item
If $f(\eta)=z$ then $\eta$ is of right codimension for every $y\in 
\mathcal{Y}_{z}$.
\item
Assuming (a), if $\pi_{\eta} :X_{\eta}\to \mathcal{Y}_{z}, \: X_{\eta}\subset 
\mathbf{P}(\mathcal{E}_{z})$ is the Gorenstein quadruple covering determined by 
$\eta$, then $X_{\eta }$ is smooth and irreducible.
\item
Assuming (a) and (b) the discriminant scheme of $\pi_{\eta}:X_{\eta}\to 
\mathcal{Y}_{z}$ is a smooth subscheme of $\mathcal{Y}_{z}$.
\end{enumerate}
Suppose $\mathbb{H}_{0}\ne \emptyset$. Consider the base change 
$\mathcal{Y}'=\mathcal{Y}\times _{Z}\mathbb{H}_{0}$ and let 
$\mathcal{E}'=\pi_{1}^{*}\mathcal{E},\: 
\mathcal{F}'=\pi_{1}^{*}\mathcal{F}$. 
Then every $\eta_{0}\in \mathbb{H}_{0}$ has a neghborhood 
$U=f^{-1}(V)\cap \mathbb{H}_{0}$, where $V$ is a Zariski open subset of 
$Z$, such that it exists a smooth quadruple covering $\mathcal{X}_{U}\to 
\mathcal{Y}'_{U}=\mathcal{Y}\times _{Z}U$ with the property that for 
every $\eta \in U$ with $f(\eta)=z$ the fiber $\mathcal{X}_{\eta}\to 
\mathcal{Y}'_{\eta}$ is equivalent to $\pi_{\eta}:X_{\eta}\to 
\mathcal{Y}_{z}$.
\end{lem}
\begin{proof}
The statement is local with respect to $Z$. According to \cite{Ha} 
Ch.III Ex.12.4 there is an invertible sheaf $\mathcal{L}$ on $Z$ such 
that $\det \mathcal{E} \cong \det \mathcal{F} \otimes 
q^{*}\mathcal{L}$. So we may assume that $Z$ is irreducible and 
$\det \mathcal{E} \cong \det \mathcal{F}$ on $\mathcal{Y}$. The proof 
then proceeds similarly to that of Lemma 2.8 of \cite{K1} and uses 
Theorem 4.4 of \cite{CE}. 
If $\mathbb{H}_{0}$ is empty there is nothing to prove. Suppose 
$\mathbb{H}_{0}\ne \emptyset$.
\begin{stepa}
Let $\mathbb{H}'$ be the set of $\eta \in \mathbb{H}$ 
for which (a) holds. We claim $\mathbb{H}'$ is Zariski open in 
$\mathbb{H}$. 
Let $\rho:\mathbf{P}(\mathcal{E})\to \mathcal{Y}$ be the projectivization and 
let $\mathcal{N}=\rho^{*}\mathcal{F}$. One has an isomorphism 
$\varPhi:q_{*}(\Check{\mathcal{F}}\otimes 
S^{2}\mathcal{E})\overset{\sim}{\lto}(q\circ \rho)_{*}
\Check{\mathcal{N}}(2)$. 
Every $\eta \in \mathbb{H}_{z}$ determines a section 
$\varPhi_z(\eta)$ of $\Check{\mathcal{N}}(2)_{z}$ with zero set 
$D_{0}(\varPhi_z(\eta))\subset \mathbf{P}(E)_{z}$.
We consider the incidence correspondence $\Gamma \subset 
\mathbf{P}(\mathcal{E})\times _{Z}\mathbb{H}$ defined as follows. 
\[
\Gamma = \{(x,\eta)|x\in D_{0}(\varPhi_z(\eta))
\; \text{where}\; x\in 
\mathbf{P}(\mathcal{E})_{y}, \eta \in \mathbb{H}_{z}, y\in 
\mathcal{Y}_{z}\}.
\]
Consider the projection $\varepsilon :\Gamma \to \mathcal{Y}\times _Z
\mathbb{H},\; \varepsilon (x,\eta)=(y,\eta)
$. An element $\eta \in \mathbb{H}_{z}$ fails to be of right codimension 
in $y\in \mathcal{Y}_{z}$ if and only if 
$\varPhi_z(\eta)\in 
H^0(\mathbf{P}(E)_{y},\Check{\mathcal{N}}(2)_{y})\cong 
H^0(\mathbb{P}^{2},\mathcal{O}_{\mathbb{P}^{2}}(2)\oplus 
\mathcal{O}_{\mathbb{P}^{2}}(2))$ determines two degenerate conics in 
$\mathbb{P}^{2}$ with a common line. This happens if and only if the 
zero set of the section $\varPhi_z(\eta)$  has dimension $\geq 1$
Equivalently $(y,\eta)\in \Sigma
\subset \mathcal{Y}\times _{Z}\mathbb{H}$ where 
$\Sigma $ is the subset of 
points for which $\dim \varepsilon ^{-1}(y,\eta)\geq 1$. Since $\Sigma$ is 
closed in $\mathcal{Y}\times _{Z}\mathbb{H}$ and  since 
$\mathcal{Y}\times _{Z}\mathbb{H}\to \mathbb{H}$ is proper, the 
projection of 
$\Sigma$ in $\mathbb{H}$ is closed.
Thus  $\mathbb{H}'$ is open in $\mathbb{H}$.
\end{stepa}
The proof of (b) and (c) is similar to that of \cite{K1} Lemma~2.8. 
Namely if $\mathcal{E}_{\mathbb{H}}$ and $\mathcal{F}_{\mathbb{H}}$ are the pull-backs via the projection 
$\mathcal{Y}\times _{Z}\mathbb{H}\to \mathcal{Y}$ one constructs a 
tautological section 
$N\in H^0(\mathcal{Y}\times 
_{Z}\mathbb{H},\Check{\mathcal{F}}_{\mathbb{H}}\otimes S^{2}
\mathcal{E}_{\mathbb{H}})$ and 
further proceeds as in Step 2 and Step 3 of the proof of Lemma 2.8 of 
\cite{K1}.
\end{proof}

\begin{block}\label{s2.38}
{\bf End of the proof of Theorem~\ref{s2.9}.} 
Given $n=2e\geq 2$
it suffices to prove the following property
for  each of the types of $(E,F)$ considered in 
(\ref{s2.12}) -- (\ref{s2.32}). 
 Unless this is 
 the type specified in 
the theorem the set of 
equivalence classes $[X\to Y]\in \mathcal{H}^{0}_{4,n}(Y)$ with
associated pair of the given type is either empty or if nonempty it is 
contained in a closed subscheme of $\mathcal{H}^{0}_{4,n}(Y)$ of 
codimension $\geq 1$. In order to prove such a statement we use 
Lemma~\ref{s2.34a} and the 
calculation of the  $\# \, \text{moduli} \, [X\to Y]$ made in (\ref{s2.12}) -- (\ref{s2.32}). Let us
analyze one case, all others being similar. Consider the possible types with 
indecomposable $F$ of degree $e$, $E\cong E_{1}\oplus E_{2}$ where 
$E_{1}$ is indecomposable of rank 2 and degree $e_{1}$ and $E_{1}$ is 
of rank 1 and degree $e_{2},\; e=e_{1}+e_{2}$. Types which satisfy 
this condition occur as subcases of Case~1 (cf. (\ref{s2.12})) and 
Case~2 (cf. (\ref{s2.15})). 
Let $\mathcal{E}(r,d)$ be the Poincar\'{e} locally free sheaf on $Y\times J$ defined in (\ref{s2.37}) and parametrizing the indecomposable locally free sheaves on $Y$ of rank $r$ and degree $d$.
We consider 
$\mathcal{E}'=p_{12}^{*}\mathcal{E}(2,e_{1})\oplus 
p_{13}^{*}\mathcal{E}(1,e_{2})$ defined on $Y\times J\times J$ and 
$\mathcal{F}' = \mathcal{E}(2,e)$ defined on $Y\times J$. The 
invertible sheaves $\det \mathcal{E}'$ and $\det \mathcal{F}'$ induce 
respectively morphisms 
\begin{align*}
&h_{1}:J\times J\to Pic^{e}Y\to J,\quad h_{2}:J\to 
Pic^{e}Y\to J, \quad \text{where}\\
&h_{1}(u_{1},u_{2})=\det 
\mathcal{E}(2,e_{1})_{u_{1}}\otimes 
\mathcal{E}(1,e_{2})_{u_{2}}\otimes \mathcal{O}_{Y}(-ey_{0}),\\ 
&h_{2}(u) = \det \mathcal{E}(2,e)_{u}\otimes 
\mathcal{O}_{Y}(-ey_{0}).
\end{align*} 
Let $Z$ be the fibre product
\begin{equation*}
\begin{diagram}
Z           & \rTo^{f_1} & J\times J \\
\dTo^{f_2}  &            &\dTo_{h_1} \\
J           &\rTo^{h_2}  & J         \\
\end{diagram}\end{equation*}
Let $\mathcal{E}=(id\times f_{1})^{*}\mathcal{E}'$ and let 
$\mathcal{F}=(id\times f_2)^{*}\mathcal{F}'$. Since $h_{1}$ and 
$h_{2}$ are smooth, surjective morphisms the fibre product $Z$ is 
smooth and by construction $\det \mathcal{E}_{z}\cong \det 
\mathcal{F}_{z}$ for every $z\in Z$. If we specify the pair $(E,F)$ to 
be of one of the types considered in Case~1 or Case~2, e.g. $E$ has 
even degree, $F$ has odd degree ecc., then in particular 
$h^0(Y,\Check{\mathcal{F}}_{z}\otimes S^2\mathcal{E}_{z}), h^0(Y,End\, 
\mathcal{E}_{z}), h^0(Y,End\, \mathcal{F}_{z})$ are independent of 
$z\in Z$. We may now apply Lemma~\ref{s2.34a} with 
$\mathcal{Y}=Y\times Z$. If there exists $[X\to Y]\in 
\mathcal{H}^{0}_{4,n}(Y)$ which has a pair $(E,F)$ of the considered 
type then $\mathbb{H}_{0}\ne \emptyset$. Let us denote by 
$\mathbb{H}_{0}^{epi}$ the union of connected components of 
$\mathbb{H}_{0}$ which correspond to $\pi :X\to Y$ with surjective 
$\pi_*:H_1(X,\mathbb{Z})\to H_1(Y,\mathbb{Z})$. For every of the open sets 
$U$ defined in Lemma~\ref{s2.34a} consider the family of
  quadruple coverings $\mathcal{X}_{U}\to Y\times U$. 
By the 
universal property of the Hurwitz space $\mathcal{H}_{4,n}(Y)$ there is an 
associated morphism $U\to \mathcal{H}_{4,n}(Y)$. These morphisms   may be 
glued. Taking the quotient by $\mathbb{C}^{\ast}$ we 
obtain a morphism $f:\mathbb{P}\mathbb{H}_{0}^{epi}\to 
\mathcal{H}^{0}_{4,n}(Y)$. Every $[X\to Y]\in 
\mathcal{H}^{0}_{4,n}(Y)$ with a pair of the given type belongs to 
$f(\mathbb{P}\mathbb{H}_{0}^{epi})$. We may now give the precise 
meaning of Formula \eqref{es2.11}: $\# \, \text{moduli} [X\to Y] = 
\dim \overline{f(\mathbb{P}\mathbb{H}_{0}^{epi})}, \# \, \text{moduli} 
(E,F) = \dim Z$, the fibers of $\mathbb{P}\mathbb{H}_{0}^{epi}\to Z$ 
have dimension $h^0(\Check{F}\otimes S^2E)-1$, the fibers of $
f: \mathbb{P}\mathbb{H}_{0}^{epi}\to \mathcal{H}^{0}_{4,n}(Y)$ have 
dimension $h^0(End\, E)-1 + h^0(End\, F)-1$. The calculations of 
(\ref{s2.12}) and (\ref{s2.15}) show that the right-hand side of 
\eqref{es2.11} is $< n$. Therefore\; $\codim \overline{
f(\mathbb{P}\mathbb{H}_{0}^{epi})}\geq 1$. In a similar manner for 
each of the types considered in Cases 1 -- 7 one constructs a smooth 
$Z$ of dimension $\# \, \text{moduli} (E,F)$ and locally free sheaves 
$\mathcal{E}$ and $\mathcal{F}$ over $Y\times Z$ and then applies 
Lemma~\ref{s2.34a}. The calculations made in 
(\ref{s2.12}) -- (\ref{s2.32}) show that all possible types except the one 
specified in the theorem yield closed subschemes of 
$\mathcal{H}^{0}_{4,n}(Y)$ of codimension $\geq 1$. Theorem~\ref{s2.9} is 
proved.
\end{block}

\begin{thm}\label{s2.42}
Let $Y$ be an elliptic curve. Let $n$ be  a pair integer $n=2e\geq 2$. 
Let $A\in Pic^{e}Y$. The Hurwitz spaces $\mathcal{H}^{0}_{4,n}(Y)$ and 
$\mathcal{H}^{0}_{4,A}(Y)$ (cf. (\ref{s2.8})) are irreducible. The 
variety $\mathcal{H}^{0}_{4,A}(Y)$ is unirational. It is rational if 
$(e,6)=1$.
\end{thm}
\begin{proof}
Let $y_{0}\in Y$ be a fixed point. The morphism 
$h:\mathcal{H}^0_{4,n}(Y)\to Pic^{e}Y$ defined by $h([X\to Y]) = \det E$ 
is surjective with fibers $\mathcal{H}^0_{4,A}(Y)$ which are isomorphic 
to each other (cf. \cite{K1} Lemma 2.5). Hence it suffices to prove 
the statements for $\mathcal{H}^{0}_{4,A}(Y)$ with 
$A=\mathcal{O}_{Y}(ey_{0})$. We have four cases according to 
Theorem~\ref{s2.9}.

\emph{Case 1. $e\not \equiv 0(mod\, 3), e\equiv 1(mod\, 2)$, i.e. 
$(e,6)=1$}. According to Atiyah's results \cite{At} up to isomorphism 
there are unique indecomposable locally free sheaves $E$ and $F$ of 
ranks 3 and 2 respectively with $\det E \cong A \cong \det F$. Let us 
apply Lemma~\ref{s2.34a} with $Z=\{\ast\}, \mathcal{Y}=Y\times Z, 
\mathcal{E} = E, \mathcal{F} = F$. The subset $\mathbb{H}_{0}\subset 
H^0(Y,\Check{F}\otimes S^2E)$ of that lemma is not empty since by 
Theorem~\ref{s2.9} every suficiently general $[X\to Y]\in 
\mathcal{H}^{0}_{4,n}(Y)$ has an associated pair of locally free 
sheaves isomorphic to $(E,F)$. Clearly $\mathbb{H}_{0}$ is invariant 
with respect to multiplication by constants in $\mathbb{C}^{\ast}$. 
Using the universal property of the Hurwitz space 
$\mathcal{H}_{4,n}(Y)$ one obtains a morphism 
$f:\mathbb{P}\mathbb{H}_{0}\to \mathcal{H}_{4,n}(Y)$. Since $
\mathbb{P}\mathbb{H}_{0}$ is irreducible its image belongs to 
$\mathcal{H}^{0}_{4,A}(Y)$. The morphism $
f:\mathbb{P}\mathbb{H}_{0}\to \mathcal{H}_{4,n}(Y)
$ is dominant by Theorem~\ref{s2.9} and injective since $h^0(End E) = 
h^0(End F) = 1$ (cf. (\ref{s2.6})). Hence $\mathcal{H}^{0}_{4,A}(Y)$ is 
irreducible and rational.

\emph{Case 2. $e \equiv 0(mod\, 3), e\equiv 1(mod\, 2)$}. Let $F$ be 
indecomposable of rank 2 with $\det F\cong A$. In order to define $Z$ 
and $\mathcal{E}$ we recall a construction due to Friedman, Morgan and 
Witten \cite{FMW}. Let $m\geq 1$, let $\mathbb{P}^{m-1}=|my_{0}|$ and 
let $|my_{0}|_{s}$ be the open subset consisting of simple divisors. 
In Theorem 2.1 (ibid) it is constructed a locally free sheaf $U(m)$ of 
rank $m$ over $Y\times |my_{0}|$ with the property that for every 
$z\in |my_{0}|_{s}, z=y_{1}+\cdots +y_{m}, y_{i}\ne y_{j}$, the 
restriction $U(m)|_{Y\times \{z\}}\cong 
\mathcal{O}_{Y}(y_{1}-y_{0})\oplus \cdots \oplus 
\mathcal{O}_{Y}(y_{m}-y_{0})$. For Case 2 we need $m=3$. We let $Z=
|my_{0}|_{s}$ and apply Lemma~\ref{s2.34a} with $\mathcal{Y}=Y\times 
Z, q=p_{2}:Y\times Z\to Z, \mathcal{E}=U(3)
\otimes p_{1}^{*}\mathcal{O}_{Y}(\frac{e}{3}y_{0})$ and 
$\mathcal{F}=p_{1}^{*}F$. As in Case 1 we obtain a dominant morphism 
$f:\mathbb{P}\mathbb{H}_{0}\to \mathcal{H}^{0}_{4,A}(Y)$. Since 
$\mathbb{P}\mathbb{H}_{0}$ is Zariski open subset in the 
projectivization of a vector bundle over $Z$ we conclude 
$\mathcal{H}^{0}_{4,A}(Y)$ is irreducible and unirational.

\emph{Case 3. $e\not \equiv 0(mod\, 3), e\equiv 0(mod\, 2)$}. This case 
is similar to the preceeding one. We let $E$ be an indecomposable 
locally free sheaf of rank 3 with $\det E \cong A$. We consider $U(2)$ 
over $Y\times |2y_{0}|$, let $Z=|2y_{0}|_{s}$ and apply 
Lemma~\ref{s2.34a} with $\mathcal{Y}=Y\times Z, \mathcal{E}=p_{1}^{*}E,\:
\mathcal{F}=U(2)\otimes p_{1}^{*}\mathcal{O}_{Y}(\frac{e}{2}y_{0})$.

\emph{Case 4. $ e\equiv 0(mod\, 6)$}. Here we consider $U(3)$ over 
$Y\times |3y_{0}|$ and $U(2)$ over $Y\times |2y_{0}|$. We let 
$Z=|3y_{0}|_{s}\times |2y_{0}|_{s}$ and applying Lemma~\ref{s2.34a} to 
$\mathcal{Y}=Y\times Z, 
\mathcal{E}=p_{12}^{*}U(3)\otimes p_{1}^{*}\mathcal{O}_{Y}(\frac{e}{3}y_{0})$ and 
$\mathcal{F}=p_{13}^{*}U(2)\otimes p_{1}^{*}\mathcal{O}_{Y}(\frac{e}{2}y_{0})$ we obtain a dominant morphism 
$f:\mathbb{P}\mathbb{H}_{0}\to \mathcal{H}^{0}_{4,A}(Y)$. This shows 
$\mathcal{H}^{0}_{4,A}(Y)$ is irreducible and unirational.
\end{proof}

\begin{rem}\label{s2.45}
We showed in the preceeding theorem that $\mathcal{H}^{0}_{4,n}(Y)$ is 
a connected component of $\mathcal{H}_{4,n}(Y)$. The Hurwitz space 
$\mathcal{H}_{4,n}(Y)$
has three 
other connected components which correspond to quadruple coverings 
$\pi :X\to Y$ such that 
$|H_1(Y,\mathbb{Z}):\pi_{*}H_1(X,\mathbb{Z})|=2$. Namely one fixes an 
\'{e}tale covering $\pi_{2}:\tilde{Y}\to Y$ of degree 2. Then every double 
covering $\pi_{1} :X\to \tilde{Y}$ branched in $n$ points which belong 
to different fibers of $\pi_{2}$ yields a simple quadruple covering 
$\pi=\pi_{2}\circ \pi_{1} :X\to Y$. One obtains in this way a 
connected component of $\mathcal{H}_{4,n}(Y)$ isomorphic to a Zariski 
open subset of $\mathcal{H}_{2,n}(\tilde{Y})$. 
\end{rem}

\begin{rem}\label{s2.44a}
Graber, Harris  and Starr proved \cite{HGS} the 
irreducibility of the space $\mathcal{H}_{d,n}^{^{S_d}}(Y)$ parameterizing 
simple coverings with monodromy group $ S_d$ for any $Y$ of positive genus 
when $n\geq 2d$. 
\end{rem}

\begin{block}\label{s2.45a} We considered so far families of quadruple 
coverings over a fixed elliptic curve $Y$. In order to treat the 
problem of unirationality of $\mathcal{A}_{3}(1,1,4)$ we need to vary 
$Y$. Simple quadruple coverings $\pi :X\to Y$ which have 3-dimensional 
Prym varieties are branched in 6 points and according to 
Theorem~\ref{s2.9} a general $[X\to Y]$ in $\mathcal{H}^{0}_{4,6}(Y)$ 
is associated to: a pair $E=M_{1}\oplus M_{2}\oplus M_{3}$ with $\deg 
M_{i}=1$; an indecomposable $F$ with $\det F\cong \det E$ and an 
element $\eta \in H^0(Y,\Check{F}\otimes S^2E)$. Let $q:\mathcal{Y}\to 
B$ be a smooth family of elliptic curves obtained from a general 
pencil of cubic curves in $\mathbb{P}^{2}$ by blowing-up the nine base 
points and discarding the singular fibers. $B$ is an open subset of 
$\mathbb{P}^{1}$. Let $\sigma :B\to \mathcal{Y}$ be a section of the 
family, $\sigma(B) = D$. One constructs by extension as in \cite{K1} 
(2.13) a rank 2 locally free sheaf $F$ on $\mathcal{Y}$ such that for 
each $b\in B$ the restriction $F_{b}=F|_{\mathcal{Y}_{b}}$ is 
indecomposable with $\det F_{b}\cong 
\mathcal{O}_{\mathcal{Y}_{b}}(3\sigma(b))$. Folowing Section~4 of \cite{FMW} 
 let $\mathcal{V}_{3}=q_{*}\mathcal{O}_{\mathcal{Y}}(3D)$. 
The projective bundle $\mathbf{P}\mathcal{V}_{3}\to B$ has fiber 
over $b\in B$ equal to $|3\sigma(b)|\subset (\mathcal{Y}_{b})^{(3)}$. 
Let us consider the locally free sheaf $U_{0}$  over 
$\mathcal{Y}\times_{B}\mathbf{P}\mathcal{V}_{3}$ defined just after 
\cite{FMW} Theorem 4.11. It has the property that if 
$z=y_{1}+y_{2}+y_{3}\in |\mathcal{O}_{\mathcal{Y}_{b}}(3\sigma(b))|$ 
is a simple divisor, then 
\(
U_{0}|_{\mathcal{Y}_{b}\times \{z\}}\cong 
\mathcal{O}_{\mathcal{Y}_{b}}(y_{1}-\sigma(b_{1}))\oplus
\mathcal{O}_{\mathcal{Y}_{b}}(y_{2}-\sigma(b_{2}))\oplus
\mathcal{O}_{\mathcal{Y}_{b}}(y_{3}-\sigma(b_{3}))
\). 
Let $Z\subset \mathbf{P}\mathcal{V}_{3}$ be the open subset consisting 
of the simple divisors in the fibers of $\mathbf{P}\mathcal{V}_{3}$, 
let $q'=p_{2}:\mathcal{Y}'=\mathcal{Y}\times_{B}Z\to Z$ be the 
projection, let $\mathcal{E}=U_{0}\otimes 
p_{1}^{*}\mathcal{O}_{\mathcal{Y}}(D)$ and let 
$\mathcal{F}=p_{1}^{*}F$. Applying Lemma~\ref{s2.34a} and letting 
$\mathbb{P}\mathbb{H}_{0}=T$ we obtain a commutative diagram
\[
\begin{diagram}
\mathcal{X}  &  \rTo & \mathcal{Y}' & \rTo  & \mathcal{Y} \\
\dTo         &       &  \dTo         &       &   \dTo      \\
T            &  \rTo &   Z           &  \rTo &   B         \\
\end{diagram}
\]
Letting $\mathcal{Y}_{T}= \mathcal{Y}'\times _{Z}T = 
\mathcal{Y}\times_{B}T$ one obtains a family of quadruple coverings 
over $T$:
\begin{equation}\label{es2.47}
\begin{diagram}
\mathcal{X}&    &\rTo^p &       &\mathcal{Y}_{T} \\
           &\rdTo &     & \ldTo  &               \\
           &      &  T  &        &               \\
\end{diagram}
\end{equation}
\end{block}

\begin{pro}\label{s2.47}
The constructed family of quadruple coverings  has 
the following properties.
\renewcommand{\theenumi}{\alph{enumi}}
\begin{enumerate}
\item
Every sufficiently general elliptic curve is isomorphic to a fiber of \linebreak
$\mathcal{Y}\to B$.
\item
Let $b\in B$. The fibers $\mathcal{X}_{\eta}\to \mathcal{Y}_{b}$ with 
$\eta \in T,\: \eta \mapsto b$ correspond to a Zariski open nonempty subset 
of the Hurwitz space $\mathcal{H}^0_{4,A}(\mathcal{Y}_{b})$ with 
$A=\mathcal{O}_{\mathcal{Y}_{b}}(3\sigma(b))$.
\item T is a rational variety of dimension 8.
\end{enumerate}
\end{pro}
\begin{proof}
The statements follow from Lemma~\ref{s2.34a} and Theorem~\ref{s2.9}.
\end{proof}

\section{The Prym mapping}\label{s3}

\begin{block}\label{s3.48}
A family of quadruple coverings of elliptic curves is given by a 
commutative triangle
\begin{equation}\label{es3.48}
\begin{diagram}
\mathcal{X}&    &\rTo^p &       &\mathcal{Y} \\
           &\rdTo_f &     & \ldTo_q  &       \\
           &      &  T  &        &           \\
\end{diagram}
\end{equation}
where $\mathcal{X},\mathcal{Y}$ and $T$ are smooth, connected, $f$ and 
$q$ are smooth, proper of relative dimension 1 with connected fibers, 
$g(\mathcal{X}_{s})=g, g(\mathcal{Y}_{s})=1$ and $p$ is finite, 
surjective (and therefore flat) of degree 4. We will work both in the 
algebraic and the complex analytic category. Assume 
$(p_{s})_{*}:H_1(\mathcal{X}_{s},\mathbb{Z})\to 
H_1(\mathcal{Y}_{s},\mathbb{Z})$ is surjective for some (and thus for 
all) $s\in T$. Define the Prym mapping using 
Proposition~\ref{s1.1} and \cite{K1} Proposition 3.14
by $\varPhi 
:T\to \mathcal{A}_{g-1}(1,\ldots ,1,4)$ where $\varPhi(s) = [Ker(Nm_{\pi 
})]$. The Prym mapping is holomorphic and  furthermore if the family is 
algebraic then $\varPhi$ is an algebraic morphism. We aim to prove the 
unirationality of $\mathcal{A}_{3}(1,1,4)$ by proving that the family 
\eqref{es2.47} yields a dominant morphism. For this we need to verify 
that the differential of the corresponding Prym mapping is generically 
surjective.
\end{block}

\begin{block}\label{s3.49}
Given a simple quadruple covering $\pi :X\to Y$ of an elliptic curve 
such that $\pi_*:H_1(X,\mathbb{Z})\to H_1(Y,\mathbb{Z})$ is surjective 
we consider as in \cite{K1} (4.2) a commutative diagram of holomorphic 
mappings
\begin{equation}\label{es3.49}
\begin{diagram}
\mathcal{X} & \rTo^p         & \mathcal{Y} \\
\dTo^f      &                &  \dTo_q     \\
N\times H   & \rTo^{\pi_{1}} & N           \\
\end{diagram}
\end{equation}
where $q:\mathcal{Y}\to N=\Delta$ is a minimal versal deformation of 
$Y$, $H\cong \Delta^{n}$. Here $\Delta$ is a unit disk. The family of 
coverings is versal for 
deformations of $\pi :X\to Y$ in the sense of Horikawa \cite{Ho2} (cf. 
\cite{K1} Proposition~4.3). The Prym mapping may be lifted to a 
holomorphic period mapping into the Siegel upper half space

\begin{diagram}
          &                         & \mathfrak{H}_{g-1}\\
          & \ruTo^{\tilde{\varPhi}} & \dTo               \\
N\times H & \rTo^{\varPhi} & \mathcal{A}_{g-1}(D)\\
\end{diagram}
where $D=(1,\ldots,1,d)$.
Given a covering of an elliptic curve there is a canonically defined 
point $q^{-}\in |\omega_{X}|^{*}=\mathbb{P}^{g-1}$ which corresponds to 
the hyperplane $H^0(X,\omega_{X})^{-}$ of differentials which have 
trace 0. Let $s_{0}\in N\times H$ be the reference point which 
corresponds to $\pi :X\to Y$. In \cite{K1} Proposition~4.16 the following criterion is
proved :

\bigskip
\noindent
\emph{Suppose $X$ is not hyperelliptic and $g(X)\geq 4$. Then $\dim 
Ker\, d\tilde{\varPhi}(s_{0})=1$ (the minimal possible dimension) if 
and only if the point $q^{-}\in |\omega_{X}|^{*}$ does not belong to 
the intersection of quadrics which contain $\phi_{K}(X)$
}

\medskip
\noindent
In particular if $g(X)=4$ the condition is $q^{-}\notin Q$ where $Q$ 
is the unique quadric which contains $\phi_{K}(X)$.
\end{block}

\begin{pro}\label{s3.51}
Let $Y$ be an elliptic curve. Then every sufficiently general simple 
quadruple covering $\pi :X\to Y$ such that $g(X)=4$ and \linebreak
$\pi_*:H_1(X,\mathbb{Z})\to H_1(Y,\mathbb{Z})$ is surjective satisfies 
the conditions of the above criterion: $X$ is not hyperelliptic and 
$q^{-}\notin Q$.
\end{pro}
\begin{proof}
That a general quadruple cover $X$ of genus 4 is not hyperelliptic is 
a particular case of \cite{K1} Proposition~4.14. We prove that 
$q^{-}\notin Q$ for a general $[X\to Y]\in \mathcal{H}^{0}_{4,6}(Y)$ 
by way of degeneration in a manner similar to the case of triple 
coverings (cf. \cite{K1} (4.19)-(4.22)).
\begin{step}
 Choose three points $\{b_{1},b_{2},b_{3}\}\subset 
Y$. Let $C_{1}=Y,\: p_{1}:C_{1}\to Y$ be the identity mapping and let 
$x_{i}=p_{i}^{-1}(b_{i})$. Let $p_{2}:C_{2}\to Y$ be a cyclic 
unramified covering of degree 3. Let 
$p_{2}^{-1}(b_{i})=\{y_{i},y'_{i},y''_{i}\}$. We consider the 
following quadruple covering of $Y$:
\[
X'\ =\ C_{1}\sqcup C_{2}/\{x_{i}\sim y_{i}\}_{i=1}^{3},
\qquad \pi'=p_{1}\cup p_{2}:X'\to Y.
\]
The same argument as that of the proof of \cite{K1} Proposition~4.20 shows 
that $\phi_{K}(X')$ is contained in a unique quadric $Q$ which is 
reducible and $q^{-}\notin Q$.
\end{step}
\begin{step}
We construct curves $X'$ of the considered type on an elliptic 
ruled surface. Let $\eta$ be a point of order 3 in $Pic^{0}Y$. Let 
$G=\mathcal{O}_{Y}\oplus \eta$. Consider the ruled surface 
$W=\mathbf{P}(\mathcal{O}_{Y}\oplus \eta),\; \rho:W\to Y$. Let $Y_{0}$ 
be the section corresponding to $\mathcal{O}_{Y}\oplus \eta\to \eta\to 
0$. We have $Y_{0}\in |\mathcal{O}_{\mathbf{P}(G)}(1)|$. Let $y\in Y$. 
Using the notation of \cite{Ha} Ch.V\S 2 consider $D=Y_{0}+yf$. Some 
useful formulas that we use below may be found in the paper 
\cite{BLi} pp.175,176. One has $h^{1}(D)=0, h^{0}(D)=2, D^{2}=2$. Let 
$\mathfrak{e}$ be a divisor on $Y$ defined by $\eta = \wedge^{2}G\cong 
L(\mathfrak{e})$. Consider the section $Y_{\infty}$ which corresponds to 
the second 
normalization  $G\otimes \eta^{-1}$ or equivalently to $G\to 
\mathcal{O}_{Y}\to 0$. Then $Y_{\infty}\sim Y_{0}-\mathfrak{e}f,\: 
Y_{0}\cdot Y_{\infty} = 0$ (cf. \cite{Ha} Ch.V Proposition~2.9) and if $y_{1}\sim y-\mathfrak{e}$ one has 
$Y_{0}+yf\sim Y_{\infty}+y_{1}f$. This shows that the pencil 
$|D|=|Y_{0}+yf|$ has two base points, $Bs|D| = \{P_{1},P_{2}\}$ where 
$P_{1}=Y_{\infty}\cap \rho^{-1}(y),\: 
P_{2}=Y_{0}\cap\rho^{-1}(y_{1})$. From Bertini's theorem it is clear 
that the general member of $|D|$ is irreducible and nonsingular.

Now, consider $|3Y_{0}|$. We have 
$\rho_{*}\mathcal{O}_{\mathbf{P}(G)}(3) = \mathcal{O}_{Y}\oplus \eta 
\oplus \eta^{2}\oplus \mathcal{O}_{Y}$. Thus 
$h^0(W,\mathcal{O}_{W}(3Y_{0}))=2$. Since $\eta \cong L(\mathfrak{e})$ 
we have $3\mathfrak{e}\sim 0$. Thus $3Y_{0}\sim 3Y_{\infty}$. We 
conclude $|3Y_{0}|$ is a pencil without base points. Let $\varphi = 
\varphi_{|3Y_{0}|}:W\to \mathbb{P}^{1}$. Let $W\to 
Z\overset{g}{\lto}\mathbb{P}^{1}$ be the Stein decomposition.  
Every fiber of $\rho:W\to \mathbb{P}^{1}$ maps surjectively to $Z$, so 
$Z\cong \mathbb{P}^{1}$. Since $3Y_{0}\sim 3Y_{\infty}$, if 
$\deg(g)>1$ then $deg(g)=3$ and $g:Z\to \mathbb{P}^{1}$ would have 
total ramification at the points corresponding to $3Y_{0}$ and 
$3Y_{\infty}$. Since $Z\cong \mathbb{P}^{1}$ this would imply 
$Y_{0}\sim Y_{\infty}$ which is absurd. We conclude $|3Y_{0}|$ is a 
pencil without base points whose general member is irreducible and 
nonsingular. The arithmetic genus $p_{a}(3Y_{0})=1$, so every 
sufficiently general curve $C\in |3Y_{0}|$ is an elliptic curve which 
is \'{e}tale triple covering of $Y$. Let $C_{1}\in |Y_{0}+yf|$ and 
$C_{2}\in |3Y_{0}|$ be suffuciently general. Then $C_{1}\cdot C_{2}=3$ 
and the intersection points belong to different fibers of $C_{2}\to 
Y$ since $C_{1}\cdot f =1$. Furtermore one may choose $C_{1}$ and 
$C_{2}$ in such a way that $C_{1}\cap C_{2}$ does not have points in 
common with $Bs|Y_{0}+yf|$. We see that $X'=C_{1}\cup C_{2}$ is a 
 curve of the type considered in Step 1. It is easy to show by 
Bertini's theorem that $X'$ is smoothable. However we need a stronger 
statement.
\end{step}
\begin{step}
\emph{$X'$ is strongly smoothable} (cf. \cite{HH} p.100). This 
means that there exist smooth, connected $\mathcal{X}$ and $B$, $\dim 
\mathcal{X} = 2,\: \dim B = 1$ and an embedding $\mathcal{X}\subset 
W\times B$ such that the second projection $\pi_{2}:\mathcal{X}\to B$ 
is proper, one of its fibers $\pi_{2}^{-1}(b_{0})\cong X'$ and all 
other fibers are smooth curves in $W$. We notice that the standard 
criterion for such smoothing $K_{W}\cdot C_{i}<0$ for every $i$ 
(cf. \cite{No}) cannot be applied here since $K_{W}\cdot C_{2} =0$. We use 
the smoothing technique  of Hartshorne and Hirschowitz \cite{HH}. In 
their paper it is stated for curves in $\mathbb{P}^{3}$, however the 
arguments can be easily extended to curves lying in arbitrary smooth 
projective variety, in particular to the simple case of curves on a 
surface. Let $N_{X'}$ be the normal sheaf of $X'=C_{1}\cup C_{2}$ and 
let $T^{1}_{X'}$ be the $T^{1}$ functor of Lichtenbaum-Schlessinger 
\cite{LS}. Since $X'$ is a curve whose singular points $P\in S=Sing\, 
X$ are nodes one has $T^{1}_{X'}\cong \oplus_{P\in S} T^{1}_{P}$ where 
$T^{1}_{P}\cong \mathbb{C}_{P}$. According to \cite{HH} 
Proposition~1.1 we have to prove two things: a) $H^{1}(N_{X'})=0$ and 
b) the natural map $H^0(N_{X'})\to H^0(T^{1}_{P})$ is surjective for 
every node $P$. According to \cite{HH} Corollary~3.2 for $i=1,2$ one 
has an exact sequence
\[
O\lto N_{C_{i}}\lto N_{X'}|_{C_{i}} \lto T^{1}_{S} \lto 0.
\]
Hence $N_{X'}|_{C_{1}}\cong N_{C_{1}}(x_{1}+x_{2}+x_{3}),\; 
N_{X'}|_{C_{2}}\cong N_{C_{2}}(y_{1}+y_{2}+y_{3})$. We have $\deg 
N_{C_{1}}=2,\; N_{C_{2}}\cong \mathcal{O}_{C_{2}}$. We now apply the 
analog of \cite{HH} Theorem~4.1. The conditions that we have to check 
in our 
case are: a) $H^{1}(C_{2},N_{C_{2}}(y_{i}+y_{j}))=0$ for $1\leq i<j\leq 
3$ and b) $H^{1}(C_{1},N_{C_{1}}) = 0$. Both are obviuosly satisfied. 
The same proof as that of \cite{HH} Theorem 4.1 then shows $X'$ is 
strongly smoothable.
\end{step}
\begin{step}
Let $p:\mathcal{X}\to Y\times B$ be the composition of 
$\mathcal{X}\hookrightarrow W\times B$ with $\rho \times id: W\times 
B\to Y\times B$. We claim that replacing $B$ by a neighborhood of 
$b_{0}$ we may assume that for each $b\in B-\{b_{0}\}$ the covering 
$p_{b}:\mathcal{X}_{b}\to Y\times \{b\}$ is simple and 
$(p_{b})_{*}: H_1(\mathcal{X}_{b},\mathbb{Z})\to H_1(Y,\mathbb{Z})$ is 
surjective. The first statement is clear since $\pi' :X'\to Y$ has the 
property that each fiber $(\pi')^{-1}(y)$ has at least 3 elements. The 
secod property is topological, so it suffices to prove it replacing 
$B$ by a small disk $\Delta$ of $b_{0}$. Let $b\in \Delta - 
\{b_{0}\}$. We have a commutative diagram
\begin{diagram}
\mathcal{X}_b &\quad \overset{i}{\hookrightarrow} 
\quad & \mathcal{X} & \quad \hookleftarrow  \quad & \mathcal{X}_{b_0} \\
\dTo          &                 &   \dTo      &                   & \dTo \\
Y\times \{b\} &\quad \hookrightarrow \quad&  Y\times \Delta & \quad \hookleftarrow \quad& Y\times \{b_0\}\\
\end{diagram} 
It is well-known that $\mathcal{X}$ is a deformation retraction of 
$\mathcal{X}_{b_{0}}$ and $i_{*}:H_1(\mathcal{X}_{b},\mathbb{Z})
\linebreak
\to 
H_1(\mathcal{X},\mathbb{Z})$ is surjective \cite{Cl,Ho}. By the
Mayer-Vietoris sequence it is clear that 
$\pi'_{*}:H_1(\mathcal{X}_{b_{0}},\mathbb{Z})\to H_1(Y,\mathbb{Z})$ is 
surjective. Thus 
$(p_{b})_{*}: H_1(\mathcal{X}_{b},\mathbb{Z})\to H_1(Y,\mathbb{Z})$ is 
surjective.
\end{step}
\begin{step}
Replacing $B$ by a neighborhood of $b_{0}$ we may 
achieve that for each $b\in B-\{b_{0}\}$ the canonical map 
$\phi_{K}:\mathcal{X}_{b}\to \mathbb{P}^{3}$ is an embedding and the 
point $q^{-}(b)$  corresponding to 
$H^0(\mathcal{X}_{b},\omega_{\mathcal{X}_{b}})^{-}$ does not belong to 
the unique quadric which contains $\phi_{K}(\mathcal{X}_{b})$. This 
statement is proved by the same argument used in \cite{K1} 
Proposition~4.22.
\end{step}
We may now conclude the proof of the proposition. Let 
$T=\mathcal{H}^{0}_{4,6}(Y)$ and let (abusing the notation) 
$p:\mathcal{X}\to Y\times T$ be the corresponding universal family of 
quadruple coverings. We showed in Steps 4 and 5 that there is an $s\in 
T$ such that the statement of the proposition holds for 
$p_{s}:\mathcal{X}_{s}\to \mathcal{Y}\times \{s\}$. Applying again to 
this family the argument of \cite{K1} Proposition~4.22 and using the 
fact that $T$ is irreducible (Theorem~\ref{s2.42}) we conclude 
that the property $q^{-}\notin Q$ holds for a Zariski open dense 
subset of $\mathcal{H}^{0}_{4,6}(Y)$.
\end{proof}

\noindent
Our approach yields an alternative proof of a result due to Birkenhake, 
Lange and van Straten \cite{BLvS}.

\begin{thm}\label{s4.56}
The moduli space of polarized abelian surfaces
$\mathcal{A}_{2}(1,4)$ is unirational.
\end{thm}
\begin{proof}
This is proved in the same way as \cite{K1} Theorem~5.1 by fixing 
$A\in Pic^{2}Y$, proving that the Prym mapping 
$\varPhi:\mathcal{H}^{0}_{4,A}(Y)\to \mathcal{A}_{2}(1,4)$ is dominant 
and thus deducing the unirationality of $\mathcal{A}_{2}(1,4)$ from the 
unirationality of $\mathcal{H}^{0}_{4,A}(Y)$ proved in 
Theorem~\ref{s2.42}.
\end{proof}

\begin{thm}\label{s3.58}
The moduli spaces of polarized abelian threefolds 
$\mathcal{A}_{3}(1,1,4)$ and $\mathcal{A}_{3}(1,4,4)$ are unirational.
\end{thm}
\begin{proof}
By a result of Birkenhake and Lange \cite{BL3} the moduli spaces 
$\mathcal{A}_{3}(1,1,4)$ and $\mathcal{A}_{3}(1,4,4)$ are isomorphic 
to each other. So it suffices to prove that $\mathcal{A}_{3}(1,1,4)$ 
is unirational. The proof is analogous to the proof of the 
unirationality of $\mathcal{A}_{3}(1,1,3)$ (cf. 
\cite{K1} Theorem~5.3). Using the family of Proposition~\ref{s2.47} 
one obtains the Prym morphism $\varPhi:T\to \mathcal{A}_{3}(1,1,4)$. 
We wish to prove $\varPhi$ is dominant. If $s_{0}\in T$ is general 
enough then using Proposition~\ref{s3.51} and the same argument as in 
\cite{K1} Theorem~5.3 one obtains 
a lifting of $\varPhi$,
 \quad 
$\tilde{\varPhi}':S\to \mathfrak{H}_{3}$
in a complex neighborhood $S$ of 
$s_{0}$
   and a commutative diagram of holomorphic mappings
\begin{diagram}
S          & \rTo^{\tilde{\varPhi}'} & \mathfrak{H}_{3} \\
\dTo^{\mu} & \ruTo_{\tilde{\varPhi}} &                  \\
N\times H  &                         &                   \\
\end{diagram}
such that a neigborhood of $\tilde{\varPhi}'(s_0)$ is contained in 
$\tilde{\varPhi}(N\times H)$. Consider the family of quadruple 
coverings $\mathcal{X}\to \mathcal{Y}\times_{N}(N\times H)$ induced 
from \eqref{es3.49}. Let $\mathcal{E}$ and $\mathcal{F}$ be the 
associated locally free sheaves of rank 3 and 2. The set of $u\in 
N\times H$ such that $\mathcal{F}_{u}$ is stable and $\mathcal{E}_{u}$ 
is semistable is open in $N\times H$ (cf. \cite{K1} Appendix~B).
Moreover 
it is open the set of $u\in 
N\times H$ such that $\mathcal{E}_{u}$ is regular polystable (ibid.)
Since $\deg 
\mathcal{E}_{u}=3=rk\, \mathcal{E}_{u}$, being regular polystable means that 
$\mathcal{E}_{u}$ is isomorphic to a direct sum $M_{1}\oplus 
M_{2}\oplus M_{3}$, where $M_i$ are invertible sheaves which are pairwise non-isomorphic to each other. So for every $u$ in a neighborhood of $\mu(s_0)$ the pair $(\mathcal{E}_{u},\mathcal{F}_{u})$ is
of the type specified in Theorem~\ref{s2.9}.
 Composing a covering of an 
elliptic curve
with its translations  one obtains coverings which have the same kernel 
of the norm map of the Jacobians while the determinants of the 
Tschirnhausen modules vary over the whole Picard group of invertible sheaves of degree $-e$, in our case 
$e=3$. By the construction of the family of Proposition~\ref{s2.47} 
this shows that $\tilde{\varPhi}'(S)$ contains the image by 
$\tilde{\varPhi}$ of a certain neighborhood of $\mu(s_{0})$ in 
$N\times H$. Therefore $\tilde{\varPhi}'(S)$ contains a neighborhood 
of $\tilde{\varPhi}'(s_{0})$ in $\mathfrak{H}_{3}$. This shows 
$\varPhi:T\to \mathcal{A}_{3}(1,1,4)$ is dominant. Therefore $
\mathcal{A}_{3}(1,1,4)$ is unirational.
\end{proof}

\bigskip
\noindent
The proof of the theorem as well as \cite{K1} Proposition~3.14 yield 
the following corollary.

\begin{cor}\label{s3.60}
Every sufficiently general abelian threefold  with polarization of 
type $(1,1,4)$ is isomorphic to the Prym variety of a simple quadruple 
covering  of an elliptic curve branched in 6 points. Every sufficiently 
general abelian threefold with polarization of type $(1,4,4)$ is 
isomorphic to $Pic^{0}X/\pi^{*}Pic^{0}Y$ for a certain qudruple covering
$\pi :X\to Y$ as above.
\end{cor}

\bibliographystyle{amsalpha}
\providecommand{\bysame}{\leavevmode\hbox to3em{\hrulefill}\thinspace}

\end{document}